# Riemann solver with internal reconstruction (RSIR) for compressible single-phase and non-equilibrium two-phase flows

Quentin Carmouze[1,2], Richard Saurel [1,3] and Emmanuel Lapebie [4]

[1] Aix Marseille Univ, CNRS, Centrale Marseille, LMA, Marseille, France
[2] University of Nice, LJAD UMR CNRS 7351, Parc Valrose, 06108 Nice Cedex, France
[3] RS2N SAS, Saint Zacharie, France
[4] CEA Gramat, France

**Abstract**

A new Riemann solver is built to address numerical resolution of complex flow models. The research direction is closely linked to a variant of the Baer and Nunziato (1986) model developed in Saurel et al. (2017a). This recent model provides a link between the Marble (1963) model for two-phase dilute suspensions and dense mixtures. As in the Marble model, Saurel et al. system is weakly hyperbolic with the same 4 characteristic waves, while the system involves 7 partial differential equations. It poses serious theoretical and practical issues to built simple and accurate flow solver. To overcome related difficulties the Riemann solver of Linde (2002) is revisited. The method is first examined in the simplified context of compressible Euler equations. Physical considerations are introduced in the solver improving robustness and accuracy of the Linde method. With these modifications the flow solver appears as accurate as the HLLC solver of Toro et al. (1994). Second the two-phase flow model is considered. A locally conservative formulation is built and validated removing issues related to non-conservative terms. However, two extra major issues appear from numerical experiments: The solution appears not self-similar and multiple contact waves appear in the dispersed phase. Building HLLC-type or any other solver appears consequently challenging. The modified Linde (2002) method is thus examined for the considered flow model. Some basic properties of the equations are used, such as shock relations of the dispersed phase and jump conditions across the contact wave. Thanks to these ingredients the new Riemann solver with internal reconstruction (RSIR), modification of the Linde method, handles stationary volume fraction discontinuities, presents low dissipation for transport waves and handles shocks and expansion waves accurately. It is validated on various test problems showing method's accuracy and versatility for complex flow models. Its capabilities are illustrated on a difficult two-phase flow instability problem, unresolved before.

# I. Introduction

The present contribution addresses building of a robust Riemann solver with limited dissipation for complex flow models. The reconstruction method of Linde (2002) is revisited and improved in terms of accuracy and robustness. This effort is mainly motivated by the numerical approximation of a two-phase non-equilibrium flow model developed by the authors that involves a series of theoretical challenges, presented in the following. In the present introduction the modelling context is recalled first, and the numerical approach is introduced secondly.

**Modelling context**

It is well accepted that hyperbolic models are mandatory to deal with phenomena involving wave propagation. This is the case for multiphase flows in many situations such as in particular shocks and detonations propagation in granular explosives and in fuel suspensions, as well as liquid-gas mixtures with bubbles, cavitating and flashing flows, as soon as motion is intense and governed by pressure gradients. This is thus the case of most unsteady two-phase flow situations.

Wave propagation is important as it carries pressure, density and velocity disturbances. Sound propagation is also very important as it determines critical (choked) flow conditions and associated mass flow rates. It has also fundamental importance on sonic conditions of detonation waves when the two-phase mixture is exothermically reacting (Petitpas et al., 2009).

Hyperbolicity is also related to the causality principle, meaning that initial and boundary conditions are responsible of time evolution of the solution. When dealing with first-order partial differential equations it means that the Riemann problem must have a solution, and the Riemann problem is correctly posed only if the equations are hyperbolic.

However, only a few two-phase flow models are hyperbolic in the whole range of parameters. The Baer and Nunziato (1986) model (BN) seemed the only formulation able to deal with such requirement. However, in the dilute limit at least, the acoustic properties of this model seemed inconsistent (Lhuillier et al., 2013). Indeed, with this model, the dispersed phase sound speed corresponds to the one of the pure phase, while this phase is not continuous and unable to propagate sound in reality, at least at a scale larger than particle's one. When the phase is not continuous (dispersed drops in a gas, dispersed bubbles in a liquid), the associated sound speed should vanish, such effect being absent in the formulation.

In the low particle's concentration limit, the Marble (1963) model is preferred. This model corresponds to the Euler equations with source terms for the gas phase and pressureless gas dynamic equations for the particle phase (see also Zeldovich, 1970). This model is thermodynamically consistent and hyperbolic as well, except that the particle phase equations are weakly hyperbolic. In this model, contrarily to the BN model, sound doesn't propagate in the particles phase, this behaviour being more physical in this limit. However, the Marble model has a limited range of validity as the volume of the dispersed phase is neglected, this assumption having sense only for low (less than per cent) condensed phase volume fraction.

Recently, the gap between these two models has been filled (Saurel et al., 2017a). Modifying the volume fraction equation in the BN model resulted in a flow model where sound propagates only in the carrier phase. The model is hyperbolic with same 4 wave speeds as Marble's one and is thermodynamically consistent in the stiff pressure relaxation limit. Moreover, in the stiff velocity relaxation limit, the Kapila et al. (2001) model, important for diffuse interface computations (Saurel and Pantano, 2018), is recovered.

Saurel et al. (2017a) model has been solved in the same reference with a Godunov type scheme based on Rusanov (1961) solver. However, as well known, this solver is quite diffusive for stationary discontinuities and linearly degenerate fields, such as in the present context, volume fraction discontinuities and contact waves. This is precisely the motivation of the present work, focused on the building of a Riemann solver with enhanced accuracy.

**Riemann solver with internal reconstruction (RSIR)**
In the quest of Riemann solver with low dissipation for this flow model several issues appear:
- The flow model, as most two-phase flow models, presents non-conservative terms;
- Numerical experiments of typical initial value problems (IVPs) achieved with the Godunov-Rusanov method exhibit non self-similar solutions. Such behavior appears as a combination of non-conservative terms, acting as a drag force (in differential form), and stiff pressure relaxation, mandatory for this specific flow model.
- Governing equations of the dispersed phase are hyperbolic degenerate, as a single eigenvalue is responsible for characteristic waves propagation. Therefore, it is impossible to determine a basis of eigenvectors and associated Riemann invariants. Moreover, the solution can be multivalued, as for the Marble's model (Saurel et al., 1994). It means that multiple volume fraction waves may be present in the solution, rendering the analysis and design of any Riemann solver intricate.

Several attempts for the building of approximate Riemann solver were done by the authors for this flow model on the basis of,
- characteristic relations for the carrier phase and jump conditions for the dispersed one,
- HLLC-type approximation based on a local conservative formulation, that will be presented later.

None of these attempts yielded efficient solver. The authors consequently move to another type of solver, based on internal reconstruction of intermediate states, computed from a simple and robust intercell state, such as Rusanov (1961) or HLL (Harten et al., 1983). This research direction has been investigated by Linde (2002), Miyoshi and Kusano (2005) and many others, mainly in the frame of magnetohydrodynamics equations that also involve many waves in the Riemann problem. The aim is to build two intermediate states instead of one. Doing so, the method should maintain stationary discontinuities and reduce numerical diffusion during simple transport.

In the present work the Linde method is revisited, and extra physics is embedded to enhance robustness and accuracy. This is done to the price of generality loss, in the sense that the method becomes model dependent, as most Riemann solvers.

The paper is organized as follows. The Linde (2002) method is recalled in Section II in its basic version and computational examples are shown with the Euler equations. It sometimes works perfectly, but oscillations appear depending on the initial conditions. This observation motivates insertion of extra physics in the solver, resulting in significant improvements, yielding robust and accurate solutions. Similar accuracy and robustness as the HLLC solver are observed.

In Section III the internal reconstruction method is extended to the two-phase flow model of interest. As before extra physics is inserted in the closure relations. The (trivial) Rankine-Hugoniot relations of the dispersed phase are used as well as jump conditions across contact wave of the dispersed phase. Thanks to these ingredients the flow solver becomes very efficient. One-dimensional computational examples are shown in the same section.

The method is then embedded in the DALPHADT unstructured meshes code. It is used to compute fingering instability occurring during explosive dispersion of particle clouds. Such instability seems misunderstood and not reproduced by existing flow models. Intensive experimental and numerical studies were done in this area recently, as for example in Rodriguez et al. (2013), McGrath et al. (2018), Osnes et al. (2018) and Xue et al. (2018). Thanks to the new model and present numerical method, the formation stage of this instability seems correctly predicted at least qualitatively. Computational examples are shown in Section IV. Conclusions are given in Section V.

## II- Riemann solver with internal reconstruction (RSIR) for the Euler equations

The Riemann solver with internal reconstruction is a modification of the Linde (2002) solver. The original Linde solver is recalled in the frame of the Euler equations and modifications are addressed next.

The Euler equations of compressible fluids consist in a system of conservation laws,

$$\frac{\partial U}{\partial t} + \frac{\partial F}{\partial x} = 0 \tag{II.1}$$

with $U = (\rho, \rho u, \rho E)^T$ and $F = (\rho u, \rho u^2 + p, (\rho E + p)u)^T$ where $\rho$ denotes the density, u the velocity and $E (= e + \frac{1}{2}u^2)$ the total energy. The pressure p is given by a convex equation of state (EOS), as a function for example of internal energy e and density $\rho$ : $p = p(\rho, e)$. The stiffened gas (SG) EOS will be used frequently in the present contribution as,

$p(\rho, e) = (\gamma - 1)\rho e - \gamma p_\infty$,

where $\gamma$ and $p_\infty$ are typical constants for a given fluid.

This system is strictly hyperbolic with wave's speeds $\lambda_1 = u$, $\lambda_2 = u - c$ et $\lambda_3 = u + c$. The sound speed is defined by $c = \sqrt{\left(\frac{\partial p}{\partial \rho}\right)_s}$ where s denotes the entropy.

The various reconstruction methods considered in the present paper are based on HLL (Harten et al., 1983) approximate solution, or its simplified version by Rusanov (1961). In the HLL solver, the extreme waves only are used, and are estimated, following Davis (1988) by,

$$S_L = \min(|u_L| - c_L, |u_R| - c_R) \text{ and } S_R = \max(|u_L| + c_L, |u_R| + c_R), \tag{II.2}$$

where subscripts L and R denote the left and right states in the initial data of the Riemann problem. The intermediate HLL state is a consequence of the Rankine-Hugoniot relations of System (1) applied across the left and right facing waves propagating at speeds $S_L$ and $S_R$ respectively:

$$U^*_{HLL} = \frac{F_R - F_L + S_L U_L - S_R U_R}{S_L - S_R} \tag{II.3}$$

From this state, the aim is now to reconstruct two intermediate states, as illustrated in Figure II.1 and linked through Relation (II.4):

$$(S_R - S_L) U^*_{HLL} = (S_R - S_M) U^*_R + (S_M - S_L) U^*_L. \tag{II.4}$$

Relation (II.4) can be expressed as,

$U^*_{HLL} = \omega_R U^*_R + \omega_L U^*_L$,

with $\omega_R = \frac{S_R - S_M}{S_R - S_L}$ and $\omega_L = \frac{S_M - S_L}{S_R - S_L}$.

The contact wave speed is given by,

$$S_M = \frac{p_R - p_L + (\rho u)_L (S_L - u_L) - (\rho u)_R (S_R - u_R)}{\rho_L (S_L - u_L) - \rho_R (S_R - u_R)},$$

Relation (II.4) involves two unknown states, $U^*_L$ and $U^*_R$. Consequently, an extra relation is needed. Linde (2002) postulate is considered first.

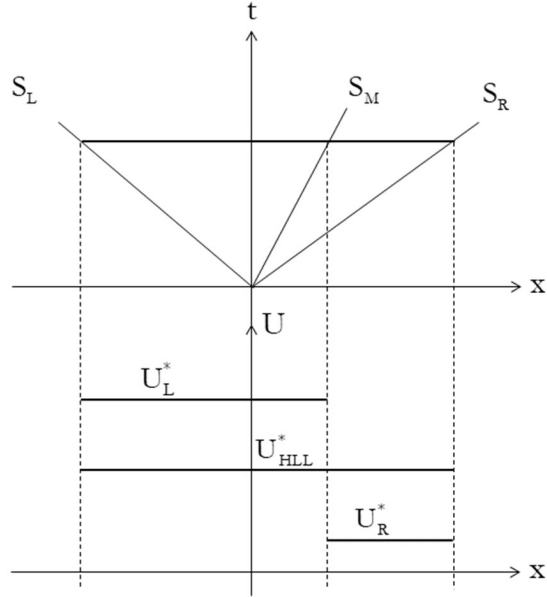

Figure II.1 – Schematic representation of the two intermediate states $U_L^*$ and $U_R^*$ rebuilt from $U_{HLL}^*$.

**II.1 Linde reconstruction**

Linde reconstruction is based on the following relation,

$$U_R^* - U_L^* = \beta(U_R - U_L) \tag{II.5}$$

where $\beta$ represents a viscosity parameter, $0 \leq \beta \leq 1$.

When $\beta$ is taken equal to zero, the HLL approximation is recovered. When $\beta = 1$ the reconstruction tends to the HLLC representation but is not equivalent, as interface conditions are ignored in this approach. It is worth to note that when $\beta = 1$, isolated density discontinuities are preserved, an important property in CFD solvers.

Relation (II.5) is then combined to Relation (II.4) resulting in:

$$\begin{cases} U_L^* = U_{HLL}^* - \omega_R \beta (U_R - U_L) \\ U_R^* = U_{HLL}^* + \omega_L \beta (U_R - U_L) \end{cases} \tag{II.6}$$

Once states $U_L^*$ and $U_R^*$ are computed the various fluxes are computed through the Rankine-Hugoniot relations,

$$\begin{cases} F_R^* = F_R + S_R (U_R^* - U_R) \\ F_L^* = F_L + S_L (U_L^* - U_L) \end{cases}. \tag{II.7}$$

Solution sampling is achieved through,

$$F^* = \begin{cases} F_L & \text{if } S_L \geq 0 \\ F_L^* & \text{if } S_L < 0 \text{ and } S_M \geq 0 \\ F_R^* & \text{if } S_R > 0 \text{ and } S_M < 0 \\ F_R & \text{if } S_R \leq 0 \end{cases}.$$

Typical solutions obtained by this solver embedded in a Godunov-type code are shown hereafter. Let us first consider transport of a density discontinuity in a uniform velocity and pressure flow with a gas governed by the ideal gas equation of state (EOS) with ɣ=1.4. Corresponding results are shown in Figure II.2.

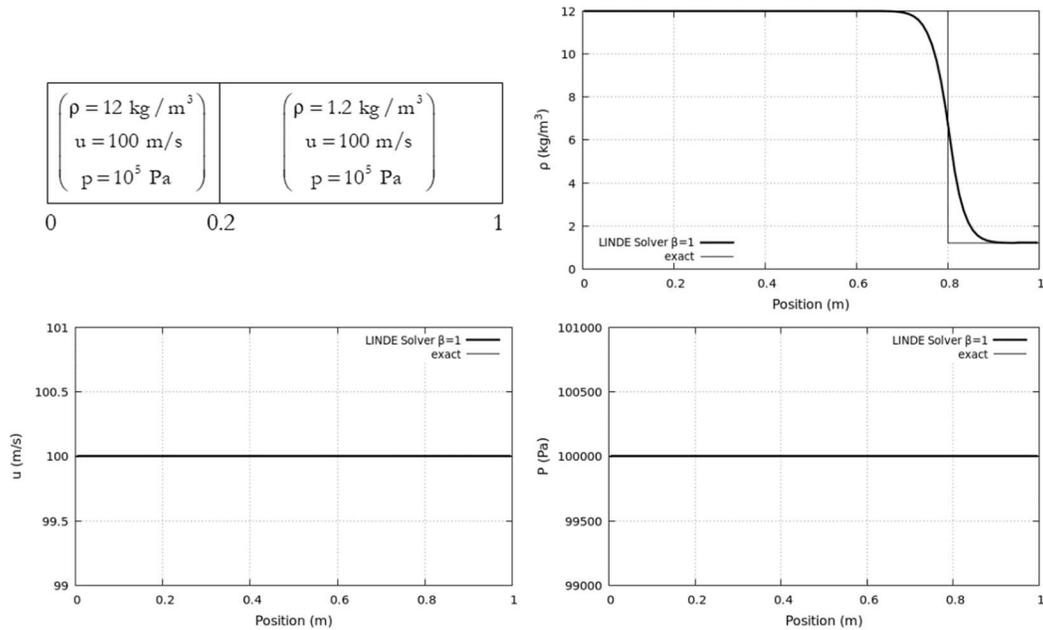

Figure II.1 – Computed results with the original Linde solver (β=1) for the transport of a density discontinuity in a uniform pressure and velocity flow. MUSCL extension of the Godunov scheme is used with Minmod limiter. 100 computational cells are used with CFL=0.5. Results are shown at time t=6ms. The density discontinuity is correctly transported, and mechanical equilibrium is maintained.

Second, a shock tube test problem is examined in Figure II.3.

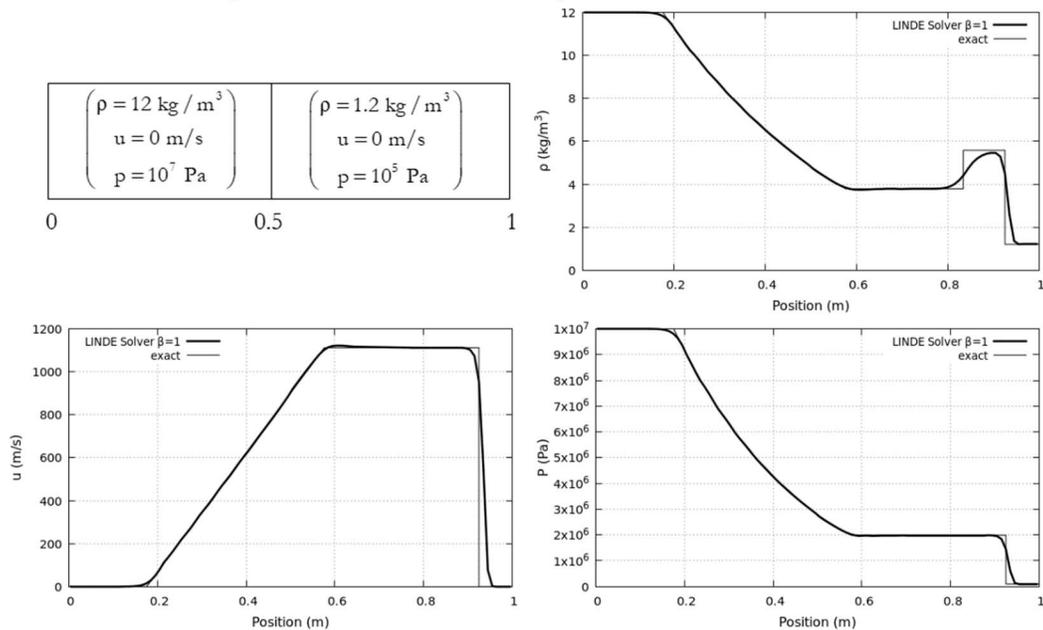

Figure II.3 – Computed results with the original Linde solver (β=1) for a shock tube test case. MUSCL extension of the Godunov scheme is used with Minmod limiter. 100 computational cells are used with CFL=0.5. Results are shown at time t=300μs. The various waves and states are computed correctly.

Third, a double expansion test is considered in Figure II.4.

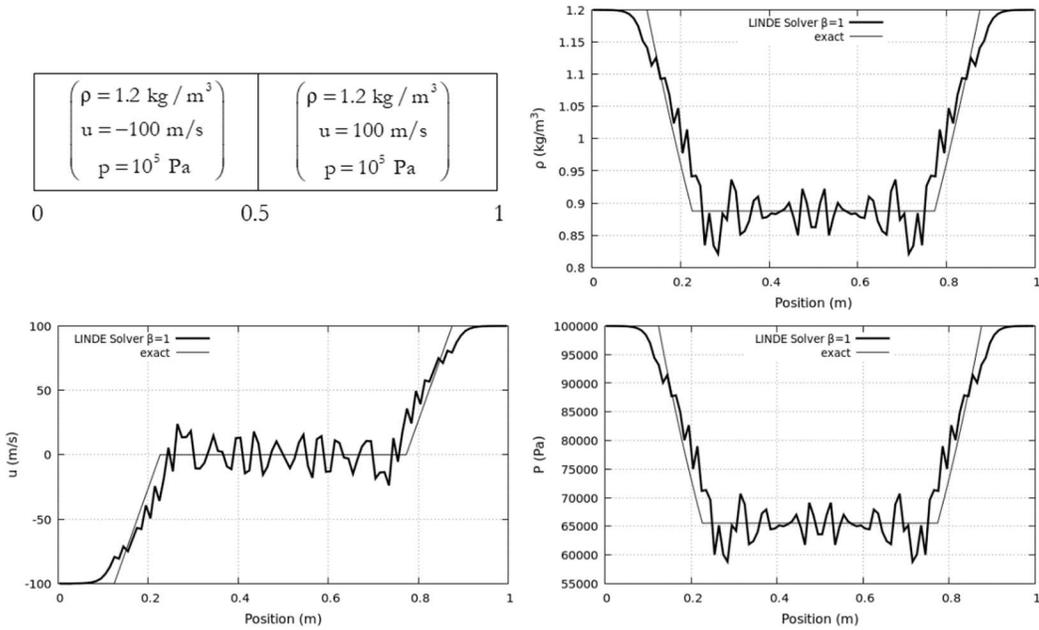

Figure II.4 – Computed results with the original Linde solver (β=1) for a double expansion test case. MUSCL extension of the Godunov scheme is used with Minmod limiter. 100 computational cells are used with CFL=0.5. Results are shown at time t=850μs. Oscillations appear. With β=0.9 same behavior of the solution is observed.

The same test is rerun in Figure II.5 with β=0.5.

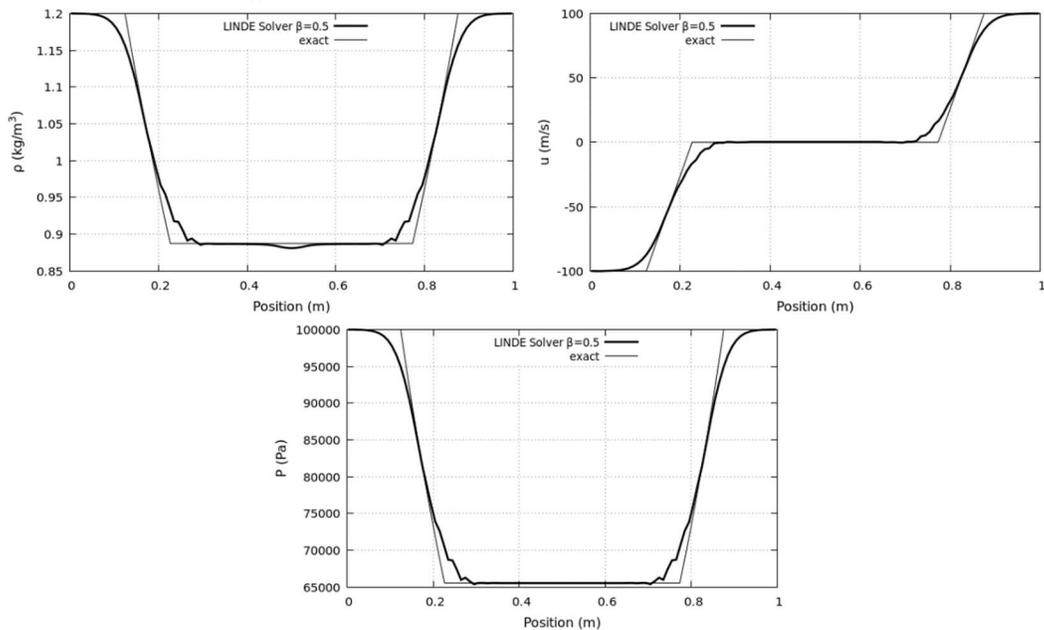

Figure II.5 - Computed results with the original Linde solver (β=0.5) for the double expansion test case of Figure II.4. MUSCL extension of the Godunov scheme is used with Minmod limiter. 100 computational cells are used with CFL=0.5. Results are shown at time t=850μs. Oscillations decrease but are still present.

Last, a double shock test is considered in Figure II.6.

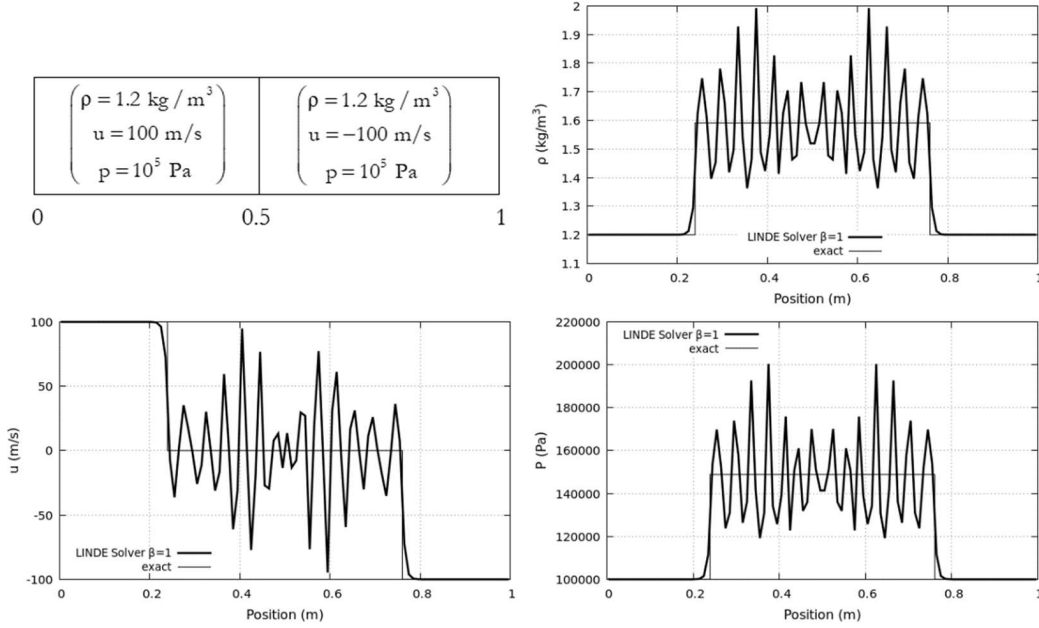

Figure II.6 – Computed results with the original Linde solver (β=1) for a double shock test case. MUSCL extension of the Godunov scheme is used with Minmod limiter. 100 computational cells are used with CFL=0.5. Results are shown at time t=850μs. Oscillations appear as before. They decrease when β tends to zero.

Accuracy of the conventional Linde method is consequently highly dependent of initial data and viscosity parameter β. It is worth to mention that in Linde (2002) a method is given to determine parameter vector for β to improve the solution. But it is important to note that as soon as β is strictly less than 1, the solver loses property of isolated stationary discontinuity preservation. Modification is thus addressed in the next paragraph.

**II.2 New reconstruction method (RSIR)**
The reconstruction derived hereafter is based on two ingredients:
- Quasi-isentropic variations across right- and left-facing waves;
- Insertion of interface conditions across the contact wave.

Although not strictly correct, in particular across strong shocks, thermodynamic evolutions through right- and left-facing waves are approximated as isentropic. Isentropic evolutions are themselves approximated through sound speed definition as,

$$c_R^2 = \frac{p_R^* - p_R}{\rho_R^* - \rho_R},$$
$$c_L^2 = \frac{p_L^* - p_L}{\rho_L^* - \rho_L}. \tag{II.8}$$

Across the contact wave $S_M$, the pressure interface condition reads,
$$p_L^* = p_R^* = p^*.$$
Thus (II.8) becomes,
$$p^* = p_R + c_R^2 \left(\rho_R^* - \rho_R\right),$$
$$p^* = p_L + c_L^2 \left(\rho_L^* - \rho_L\right).$$
Taking the difference of these two relations the following one is obtained,

$$\rho_R^* - \frac{c_L^2}{c_R^2}\rho_L^* = \left(\rho_R - \frac{c_L^2}{c_R^2}\rho_L\right) + \frac{p_L - p_R}{c_R^2}.$$

For simplicity reasons another approximation is now done with respect to the square sound speeds, that are assumed closed to a certain average:

$$c_L^2 \approx c_R^2 \approx \overline{c}^2.$$

Consequently,

$$\rho_R^* - \rho_L^* = \rho_R - \rho_L + \frac{p_L - p_R}{\overline{c}^2} \text{ with } \overline{c}^2 = \frac{c_R^2 + c_L^2}{2}. \tag{II.9}$$

This relation corresponds to a modification of the first relation of System (II.5).

To maintain flexibility of the reconstruction method, parameter β is reintroduced as,

$$\rho_R^* = \rho_L^* + \beta\left(\rho_R - \rho_L + \frac{p_L - p_R}{\overline{c}^2}\right). \tag{II.10}$$

Indeed, parameter β is convenient to control numerical viscosity.

From (II.10) similar relation is deduced for momentum and energy jumps across the contact wave. As $u_L^* = u_R^* = S_M$, (II.10) implies the following relation,

$$(\rho u)_R^* = (\rho u)_L^* + \beta\left(\rho_R - \rho_L + \frac{p_L - p_R}{\overline{c}^2}\right)S_M \tag{II.11}$$

Then, for the energy jump, using once more the pressure interface condition,

$$p_L^* = p_R^* = p^*,$$

it appears, provided that the fluid is governed by ideal gas or SG EOS,

$$(\rho e)_L^* = (\rho e)_R^*.$$

Consequently,

$$(\rho E)_R^* = (\rho e)_R^* + \left(\frac{1}{2}\rho u^2\right)_R^* = (\rho e)_L^* + \frac{1}{2}u_R^{*2}\left[\rho_L^* + \beta\left(\rho_R - \rho_L + \frac{p_L - p_R}{\overline{c}^2}\right)\right].$$

Alternatively, it reads,

$$(\rho E)_R^* = (\rho E)_L^* + \beta\left(\rho_R - \rho_L + \frac{p_L - p_R}{\overline{c}^2}\right)\frac{S_M^2}{2}. \tag{II.12}$$

Relations (II.10-II.12) summarize as,

$$U_{R,k}^* = U_{L,k}^* + \psi\Lambda_k, \text{ k=1,..,3}. \tag{II.13}$$

with,

$$\psi = \beta\left[\rho_R - \rho_L + \frac{p_L - p_R}{\overline{c}^2}\right] \text{ and } \Lambda = \left(1, S_M, \frac{S_M^2}{2}\right)^T.$$

Combining (II.4) and (II.13) the various intermediate states are rebuilt as,

$$\begin{cases} U_{L,k}^* = U_{HLL,k}^* - \omega_R \psi \Lambda_k \\ U_{R,k}^* = U_{HLL,k}^* + \omega_L \psi \Lambda_k \end{cases}, \text{ k=1,..,3}. \tag{II.14}$$

From these intermediate states fluxes are computed with the help of (II.7) as before.

The Riemann solver thus consists in (II.3) to compute the HLL state and (II.14) to rebuild the two intermediate states. It is thus particularly simple and contains parameter β to control numerical diffusion. Comparison with the HLLC solver is now addressed on various test problems. In all tests that follow, β is taken equal to 1 to minimize numerical diffusion.

The first test consists in the computation of a density discontinuity at rest, as shown in Figure II.7.

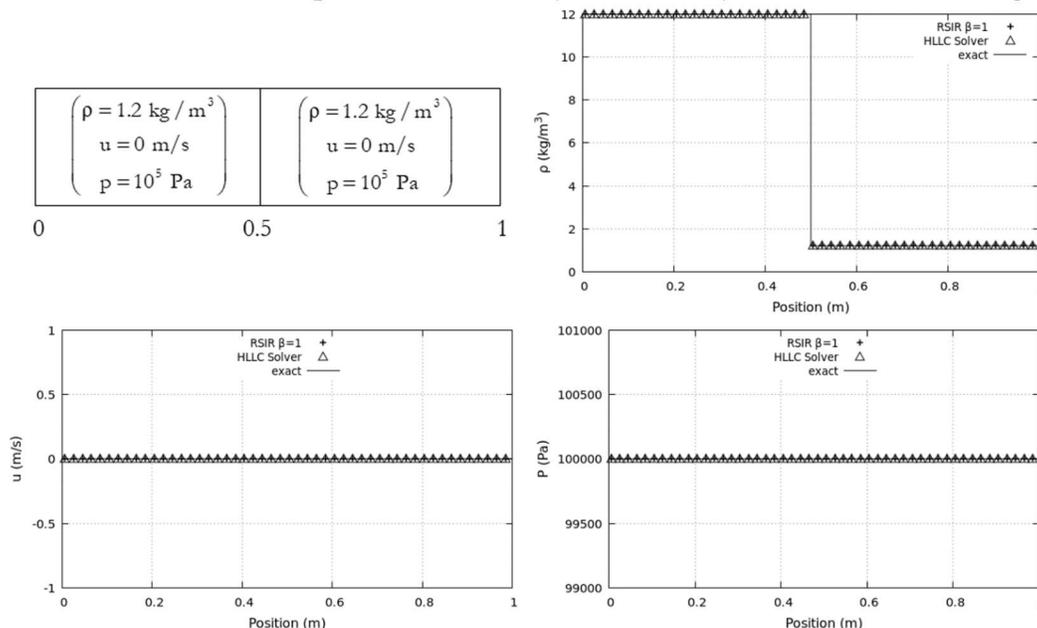

Figure II.7 – Comparison of HLLC and new solver (β=1) embedded in the Godunov-MUSCL method with Minmod limiter for the computation of a contact discontinuity at rest. 100 computational cells are used with CFL=0.5. Results are shown at time t=6ms. For the sake of clarity only 50 symbols out of 100 are plotted for both computations. Both computed results are merged.

The double-expansion and double-shock tests of Figures II.4 and II.6 are rerun with the new method and HLLC. Results are shown in Figures II.8 and II.9 showing oscillation free solutions.

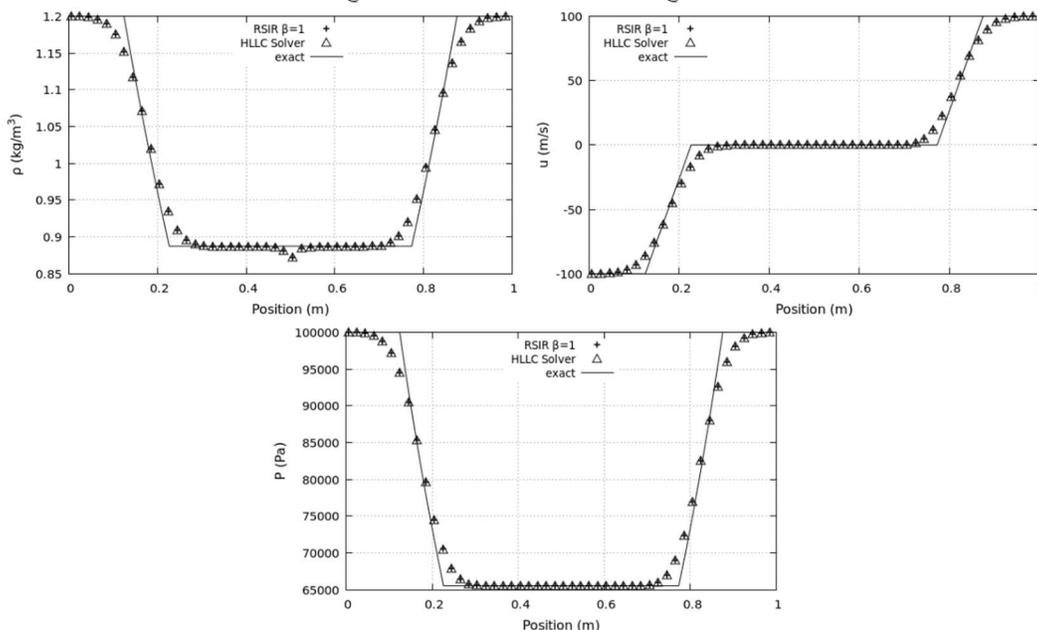

Figure II.8 – Comparison of HLLC and new solver (β=1) embedded in the Godunov-MUSCL method with Minmod limiter for the computation of a double-expansion test. 100 computational cells are used with CFL=0.5. Results are shown at time t=850μs. For the sake of clarity only 50 symbols out of 100 are plotted for both computations. Both computed results are merged.

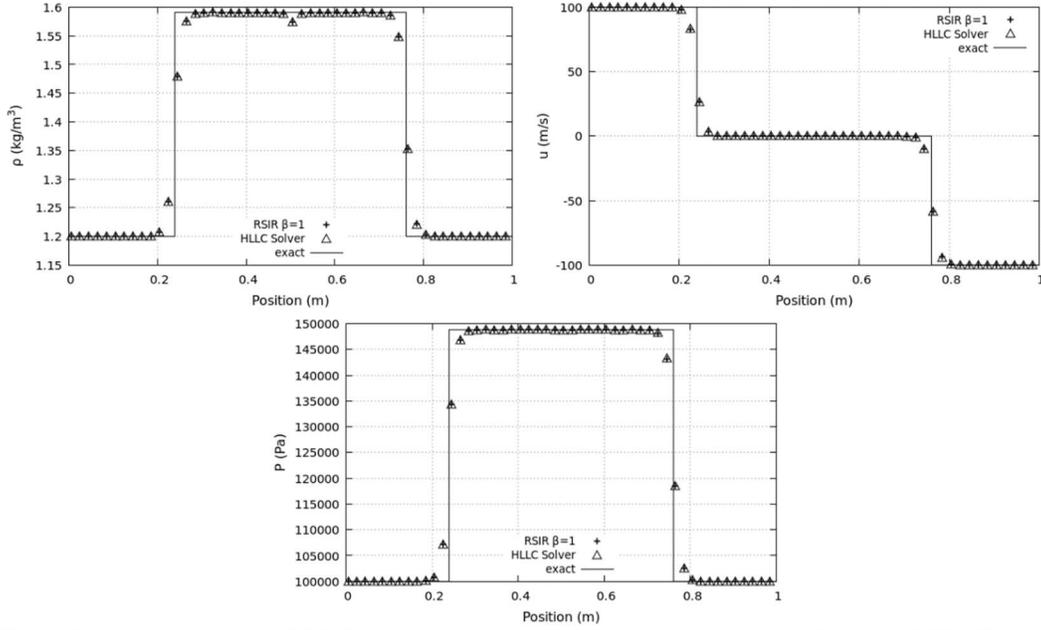

Figure II.9 – Comparison between HLLC and the new solver (β=1) embedded in the Godunov-MUSCL method with Minmod limiter for the computation of a double-shock test. 100 computational cells are used with CFL=0.5. Results are shown at time t=850µs. For the sake of clarity only 50 symbols out of 100 are plotted for both computations. Both computed results are merged.

The RSIR method is consequently efficient, at least provided that the fluid is governed by the ideal gas or SG EOS, as considered in the derivation of (II.12). Considering a more sophisticated but convex EOS, such as for example NASG (Le Metayer and Saurel, 2016),

$$p = p(\rho, e) = \frac{(\gamma - 1)\rho e}{1 - \rho b} - \gamma p_\infty,$$

where b denotes the covolume, the present method requires modification.

Combining (II.4) and (II.10) star densities $\rho_R^*$ and $\rho_L^*$ are determined. With the help of (II.8) the star pressures $p_R^*$ and $p_L^*$ are estimated. The associated internal energies are computed with the help of the EOS, $e_k^*\left(p_k^*, \rho_k^*\right)$, k=L, R.

Relation (II.12) is thus replaced by,

$$\left(\rho E\right)_R^* - \left(\rho E\right)_L^* = \rho_R^* e\left(p_R^*, \rho_R^*\right) - \rho_L^* e\left(p_L^*, \rho_L^*\right) + \left(\rho_R^* - \rho_L^*\right)\frac{1}{2}S_M^2 \qquad (II.15)$$

The RSIR solver now reads,

$$\begin{cases} U_{L,k}^* = U_{HLL,k}^* - \omega_R \psi_k \\ U_{R,k}^* = U_{HLL,k}^* + \omega_L \psi_k \end{cases}, \; k=1,..,3, \qquad (II.16)$$

with,

$$\psi = \begin{pmatrix} \beta\left(\rho_R - \rho_L + \dfrac{p_L - p_R}{\bar{c}^2}\right) \\ \beta\left(\rho_R - \rho_L + \dfrac{p_L - p_R}{\bar{c}^2}\right) S_M \\ \rho_R^* e(p_R^*, \rho_R^*) - \rho_L^* e(p_L^*, \rho_L^*) + (\rho_R^* - \rho_L^*)\dfrac{1}{2} S_M^2 \end{pmatrix}.$$

This RSIR method is compared to computations based on HLLC solver for the transport of a density discontinuity in Figure II.10 and on shock tube computations in Figure II.11. NASG EOS parameters used in these computations are those of liquid water: ɣ=4.4, p∞=6.10⁸ Pa and b=5.10⁻⁵ m³/kg.

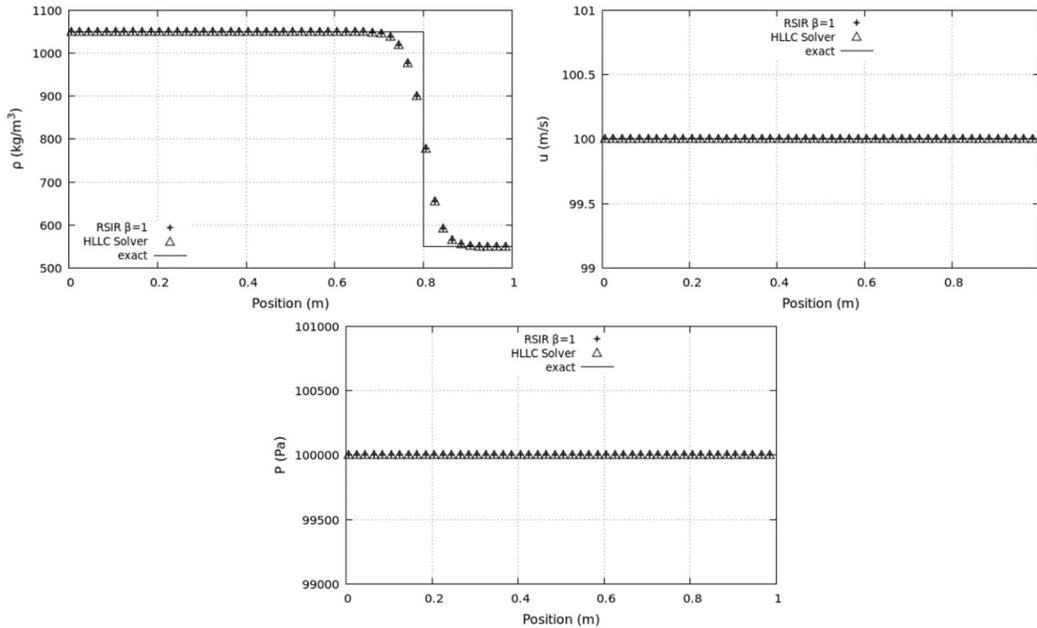

Figure II.10 – Comparison between HLLC and the new solver (β=1) embedded in the Godunov-MUSCL method with Minmod limiter for the transport of a density discontinuity in a uniform pressure and velocity flow using the NASG EOS. The initial discontinuity is located at x=0.2m. 100 computational cells are used with CFL=0.5. Results are shown at time t=6ms. For the sake of clarity only 50 symbols out of 100 are plotted for both computations. Both computed results are merged.

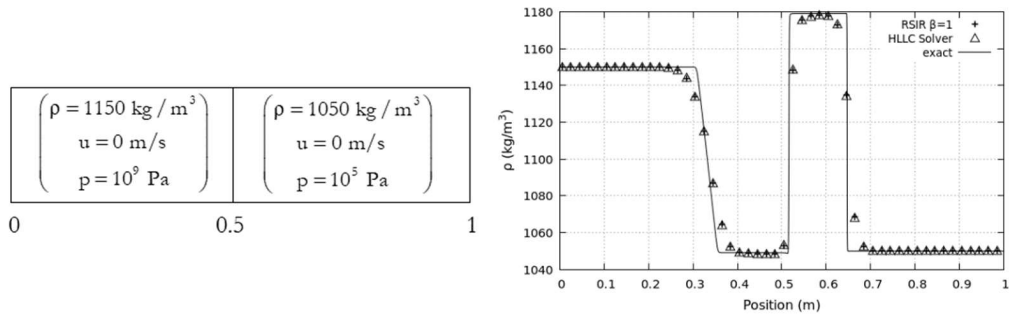

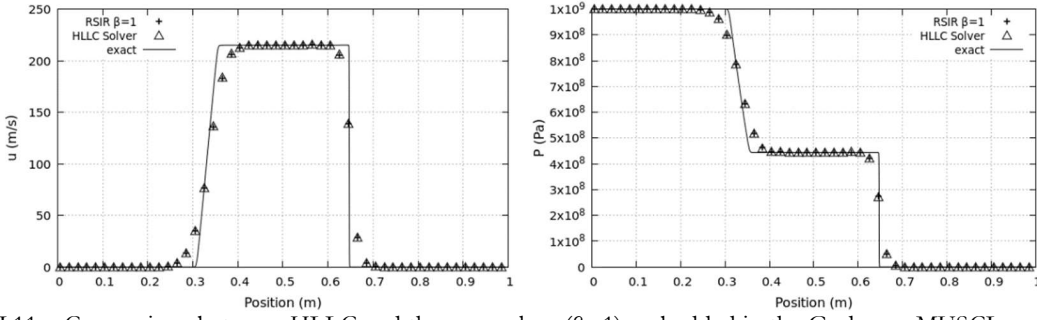

Figure II.11 – Comparison between HLLC and the new solver (β=1) embedded in the Godunov-MUSCL method with Minmod limiter for the computation of a shock tube test with liquid water and NASG EOS. 100 computational cells are used with CFL=0.5. Results are shown at time t=75µs. For the sake of clarity only 50 symbols out of 100 are plotted for both computations. Both computed results are merged.

The new solver seems as accurate as the HLLC solver, at least for the Euler equations. But it is not as flexible, or systematic, as modifications are needed when modifications of the flow model are done. Its extension to the more sophisticated flow model of Saurel et al. (2017a) is now examined. As already mentioned, this system is weakly hyperbolic. Moreover, it will be shown that solutions are not self-similar. For such a complex flow model, present authors tried but failed in the derivation of HLLC-type Riemann solver.

### III – Extension to dense-dilute two-phase flow model

The two-phase flow model of Saurel et al. (2017a), considering a dispersed phase 1 in a carrier fluid 2 is recalled hereafter. Pressure and velocity relaxation terms only are considered as interaction effects:

$$\frac{\partial \alpha_1}{\partial t} + \frac{\partial (\alpha u)_1}{\partial x} = \mu (p_1 - p_2), \qquad \mu \to +\infty$$

$$\frac{\partial (\alpha \rho)_1}{\partial t} + \frac{\partial (\alpha \rho u)_1}{\partial x} = 0,$$

$$\frac{\partial (\alpha \rho u)_1}{\partial t} + \frac{\partial (\alpha \rho u^2 + \alpha p)_1}{\partial x} = p_I \frac{\partial \alpha_1}{\partial x} + \lambda (u_2 - u_1),$$

$$\frac{\partial (\alpha \rho E)_1}{\partial t} + \frac{\partial (\alpha (\rho E + p) u)_1}{\partial x} = p_I \frac{\partial (\alpha u)_1}{\partial x} + \lambda u_I (u_2 - u_1) - \mu p_I (p_1 - p_2), \qquad (III.1)$$

$$\frac{\partial (\alpha \rho)_2}{\partial t} + \frac{\partial (\alpha \rho u)_2}{\partial x} = 0,$$

$$\frac{\partial (\alpha \rho u)_2}{\partial t} + \frac{\partial (\alpha \rho u^2 + \alpha p)_2}{\partial x} = p_I \frac{\partial \alpha_2}{\partial x} - \lambda (u_2 - u_1),$$

$$\frac{\partial (\alpha \rho E)_2}{\partial t} + \frac{\partial (\alpha (\rho E + p) u)_2}{\partial x} = -p_I \frac{\partial (\alpha u)_1}{\partial x} - \lambda u_I (u_2 - u_1) + \mu p_I (p_1 - p_2).$$

Same notations as for the Euler equations are used. Additionally $\alpha_k$ represent the volume fraction of phase k, index I is related to interfacial variables and relaxation parameters are denoted by $\lambda$ and $\mu$, with respect to velocity and pressure relaxation. Appropriate relations are given for example in Saurel et al. (2017a). It is important to mention that this flow model is valid only in the stiff pressure relaxation limit ($\mu \to +\infty$). Appropriate pressure relaxation solvers are given for example in Lallemand and Saurel (2000). This is not equivalent to strict pressure equilibrium models that are non-hyperbolic, or conditionally hyperbolic. Also, this is not a restrictive assumption for most two-phase flow applications, except possibly extreme situations, such as hot spots ignition in condensed energetic materials (Saurel et al., 2017b), where pressure relaxation is the driving effect for hot spot appearance. The present model is hyperbolic with wave speeds,

$\lambda_{1-4} = u_1$, $\lambda_5 = u_2$, $\lambda_6 = u_2 - c_2$ and $\lambda_7 = u_2 + c_2$.

These waves speeds are the same as Marble (1963) model but are significantly different to those of the Baer and Nunziato (1986) model. Combination of the equations of System (III.1) result in the following mixture entropy equation, that guarantees non-negative evolutions,

$$\frac{\partial \alpha_1 \rho_1 s_1 + \alpha_2 \rho_2 s_2}{\partial t} + \frac{\partial \alpha_1 \rho_1 s_1 u_1 + \alpha_2 \rho_2 s_2 u_2}{\partial x} = \left( \frac{(u_I - u_1)}{T_1} - \frac{(u_I - u_2)}{T_2} \right) \lambda (u_2 - u_1). \tag{III.2}$$

Indeed, admissible estimates for the interfacial velocity are,

$u_I = u_1$ or $u_I = u_2$,

or combinations of them.

System (III.1) is obviously non-conservative, as $p_I \frac{\partial \alpha_1}{\partial x}$ and $p_I \frac{\partial (\alpha u)_1}{\partial x}$ terms are present in the right-hand side of the momentum and energy equations. However, assuming $p_I = p_1$ the following Rankine-Hugoniot system is obtained (Saurel et al., 2017a):

$\alpha_1 = \alpha_1^0$,
$\rho_1 = \rho_1^0$,
$e_1 = e_1^0$,
$u_1 = u^0$, \hfill (III.3)
$\rho_2(u_2 - \sigma) = \rho_2^0(u^0 - \sigma)$,
$\rho_2 u_2(u_2 - \sigma) + p_2 = \rho_2^0 u^0(u^0 - \sigma) + p^0$,
$\rho_2 E_2(u_2 - \sigma) + u_2 p_2 = \rho_2^0 E_2^0(u^0 - \sigma) + u^0 p^0$.

These relations will be used in the RSIR derivation.

Numerical resolution of System (III.1) has been done in Saurel et al. (2017a) with the help of Rusanov (1961) solver. This solver being quite diffusive, the aim is now to build an improved solver. System (III.1) involves however three main difficulties:
- It is non-conservative;
- The eigenvalue $u_1$ is multiple. Therefore System (III.1) admits multivalued solutions (Forestier and Le Floch, 1992, Saurel et al., 1994, Bouchut et al., 2003).
- Solutions are not self-similar, as will be shown later with the help of numerical experiments.

These issues are addressed gradually in the following.

### III.1 Local conservative formulation and Rusanov-type solvers

In Saurel et al. (2017a) a Rusanov-type method was derived to determine qualitative solutions of the new flow model and validations against both exact solutions and experimental data. This method is recalled hereafter, and an improved version based on a local conservative formulation is built. The aim is to show that non-conservative terms are treated correctly through the local conservative formulation.

a) Basic Rusanov version

System (III.1) is considered in non-conservative form and in the absence of relaxation terms as,

$$\frac{\partial U}{\partial t} + \frac{\partial F(U)}{\partial x} + H\left(U, \frac{\partial U}{\partial x}\right) = 0, \tag{III.4}$$

with,

$$U = \left(\alpha_1, \ (\alpha\rho)_1, \ (\alpha\rho u)_1, \ (\alpha\rho E)_1, \ (\alpha\rho)_2, \ (\alpha\rho u)_2, \ (\alpha\rho E)_2\right)^T$$

$$F = \left((\alpha u)_1, \ (\alpha\rho u)_1, \ (\alpha\rho u^2 + \alpha p)_1, \ (\alpha(\rho E + p)u)_1, \ (\alpha\rho u)_2, \ (\alpha\rho u^2 + \alpha p)_2, \ (\alpha(\rho E + p)u)_2\right)^T$$

$$H = \left(0, \ 0, \ -p_I \frac{\partial \alpha_1}{\partial x}, \ -p_I \frac{\partial (\alpha u)_1}{\partial x}, \ 0, \ p_I \frac{\partial \alpha_1}{\partial x}, \ p_I \frac{\partial (\alpha u)_1}{\partial x}\right)^T$$

Let us denote by,

$$S = \text{Max}_k \left(|\lambda_k|_L, |\lambda_k|_R\right),$$

the maximum wave speed separating two states L and R.
The Rusanov flux reads,

$$F^* = \frac{1}{2}\left[F_R + F_L - S(U_R - U_L)\right],$$

and the Godunov scheme associated to system (III.4) reads,

$$U_i^{n+1} = U_i^n - \frac{\Delta t}{\Delta x}\left(F^*_{i+\frac{1}{2}} - F^*_{i-\frac{1}{2}}\right) + \Delta t H_i,$$

where $H_i$ is an approximation of non-conservative terms.
Following Saurel et al. (2017a), based on Saurel and Abgrall (1999) method for a slightly different flow model, approximation of these terms read;
- For the momemtum equation,

$$H_{i,u} = p_i^n \frac{\alpha^*_{1,i+\frac{1}{2}} - \alpha^*_{1,i-\frac{1}{2}}}{\Delta x} \text{ with } \alpha^*_{1,i+\frac{1}{2}} = \frac{\alpha_{1,i} + \alpha_{1,i+1}}{2}.$$

- For the energy equations,

$$H_{i,E} = p_i^n \frac{(\alpha u)^*_{1,i+\frac{1}{2}} - (\alpha u)^*_{1,i-\frac{1}{2}}}{\Delta x} \text{ with } (\alpha u)^*_{1,i+\frac{1}{2}} = \frac{(\alpha u)_{1,i+1} + (\alpha u)_{1,i}}{2}.$$

These formulas are built to respect mechanical equilibrium condition. Another version is examined hereafter.

b) Local conservative formulation

The interfacial pressure $p_I$ appears in the presence of non-conservative terms, such as $p_I \frac{\partial \alpha_k}{\partial x}$. As $p_I$ has been assumed equal to the dispersed phase pressure, the following assumption is made:

$$p_I = \begin{cases} p_{1,L} & \text{if } \alpha_{1,L} > \alpha_{1,R} \\ p_{1,R} & \text{if } \alpha_{1,L} < \alpha_{1,R} \end{cases}. \tag{III.5}$$

It means that $p_I$ is taken equal to the pressure of phase 1 when this phase is present. Possible situations are schematized in Figure III.1.

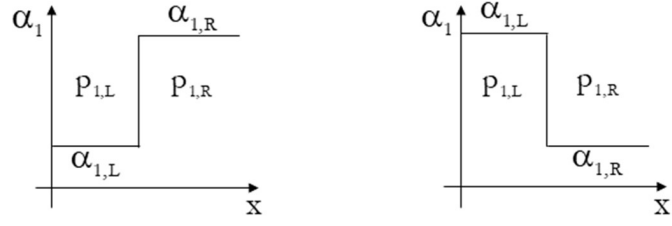

Figure III.1 – Schematic representation of the estimate for $p_I$. As a consequence of jump conditions (III.3) both volume fraction and pressure of phase 1 are invariant across the extreme waves $u_2 - c_2$ and $u_2 + c_2$ in the Riemann problem solution. Therefore, $p_I$ is taken equal to the pressure of phase 1 when this phase is present, as summarized in (III.5). As it is constant during time evolution in a given Riemann problem, it becomes a local constant.

As a consequence of Rankine-Hugoniot relations (III.3), assuming $p_I = p_1$ implies that $p_I$ becomes a local constant, as $p_1$ is invariant across right- and left-facing waves.

Thanks to this local constant, System (III.1) becomes locally conservative:

$$\frac{\partial \alpha_1}{\partial t} + \frac{\partial (\alpha u)_1}{\partial x} = 0,$$

$$\frac{\partial (\alpha \rho)_1}{\partial t} + \frac{\partial (\alpha \rho u)_1}{\partial x} = 0,$$

$$\frac{\partial (\alpha \rho u)_1}{\partial t} + \frac{\partial \left(\alpha \rho u^2 + \alpha(p - p_I)\right)_1}{\partial x} = 0,$$

$$\frac{\partial (\alpha \rho E)_1}{\partial t} + \frac{\partial \left(\alpha(\rho E + p - p_I)u\right)_1}{\partial x} = 0,$$

$$\frac{\partial \alpha_2}{\partial t} - \frac{\partial (\alpha_1 u_1)}{\partial x} = 0,$$

$$\frac{\partial (\alpha \rho)_2}{\partial t} + \frac{\partial (\alpha \rho u)_2}{\partial x} = 0,$$

$$\frac{\partial (\alpha \rho u)_2}{\partial t} + \frac{\partial \left(\alpha \rho u^2 + \alpha(p - p_I)\right)_2}{\partial x} = 0,$$

$$\frac{\partial (\alpha \rho E)_2}{\partial t} + \frac{\partial \left(\alpha(\rho E + p)u\right)_2 + (\alpha u)_1 p_I}{\partial x} = 0.$$

(III.6)

In compact form it reads,

$$\frac{\partial U}{\partial t} + \frac{\partial \Phi(U)}{\partial x} = 0,$$

with obvious definition for $\Phi(U)$.

The associated Rusanov flux is immediate,

$$\Phi^* = \frac{1}{2}\left[\Phi_R + \Phi_L - S(U_R - U_L)\right].$$

(III.7)

From $\Phi^*$, the $F^*$ flux of formulation (III.4) is deduced as,

$$F_k^* = \Phi_k^* + p_I \begin{pmatrix} 0 \\ 0 \\ \alpha_k^* \\ \Phi^*(\alpha_k) \end{pmatrix}, \quad k=1,2. \tag{III.8}$$

$\alpha_k^*$ are determined form Rusanov state as,

$$\alpha_k^* = \frac{1}{2}\left[\alpha_{k,R} + \alpha_{k,L} - \frac{\left(\Phi(\alpha_{k,R}) - \Phi(\alpha_{k,L})\right)}{S}\right]. \tag{III.9}$$

The fluxes are inserted in the same Godunov scheme as before,

$$U_i^{n+1} = U_i^n - \frac{\Delta t}{\Delta x}\left(F_{i+\frac{1}{2}}^* - F_{i-\frac{1}{2}}^*\right) + \Delta t H_i,$$

except that $H_i$ are now given by,

$$H_{i,u} = p_i^n \frac{\alpha_{1,i+\frac{1}{2}}^* - \alpha_{1,i-\frac{1}{2}}^*}{\Delta x}$$

$$H_{i,E} = p_i^n \frac{\Phi_{i+\frac{1}{2}}^*(\alpha_1) - \Phi_{i-\frac{1}{2}}^*(\alpha_1)}{\Delta x}.$$

In these expressions $\alpha_{k,i\pm\frac{1}{2}}^*$ and $\Phi_{i+\frac{1}{2}}^*(\alpha_1)$ are given by (III.9) and (III.7) respectively.

The overall scheme is consequently quite different of the basic version presented before.

The two Rusanov solvers are now compared on an arbitrary shock tube test problem, in the absence of velocity relaxation, but with stiff pressure relaxation. Corresponding results are shown in Figure III.2.

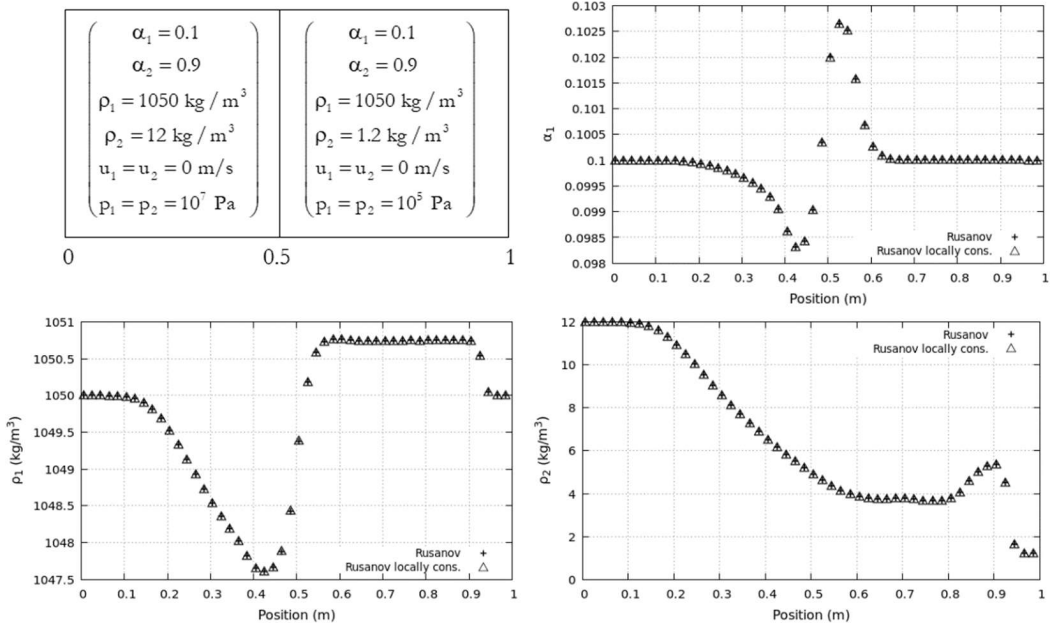

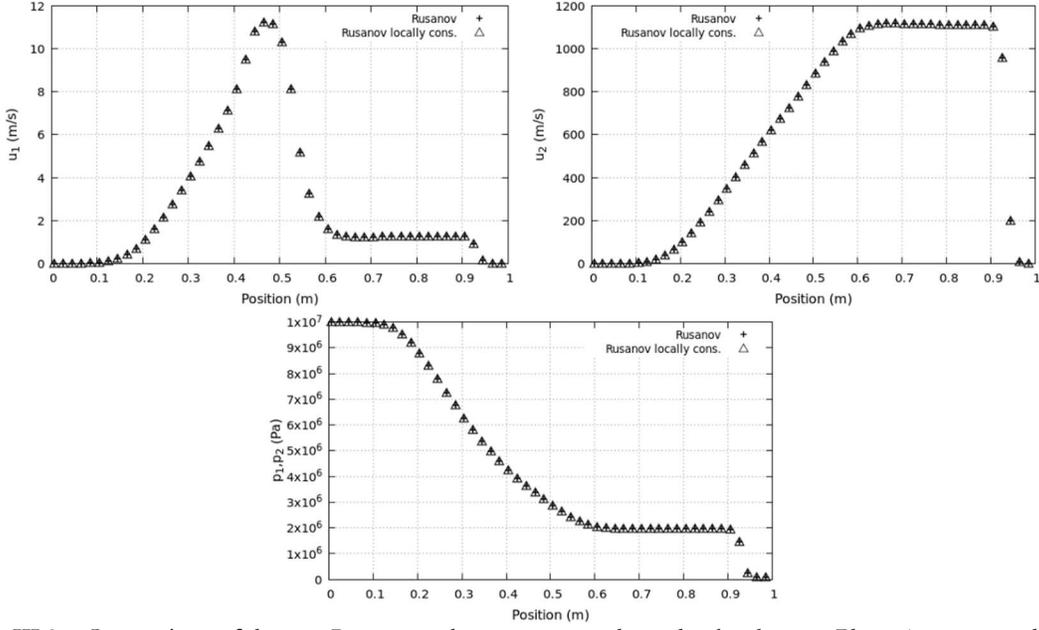

Figure III.2 – Comparison of the two Rusanov solvers on a two-phase shock tube test. Phase 1 corresponds to the dispersed fluid, considered here as liquid water, with SG EOS parameters ($\gamma_1 = 4.4$ and $p_{\infty 1} = 6\,10^8\,Pa$). Phase 2 represents the carrier phase, here air considered as ideal gas ($\gamma_2 = 1.4$). Stiff pressure relaxation is used at any time. 100 computational cells are used with CFL=0.5. Computations are done with the MUSCL reconstruction method and Minmod limiter. Results are shown at t=300µs. For the sake of clarity only 50 symbols out of 100 are plotted for both computations. Both results are perfectly merged validating the approach based on local conservative formulation.

Extra tests have been done, such as double expansion and double shock tests, always showing the same agreement. The local conservative formulation (III.6) with local constant (III.5) is consequently robust enough to be considered with the reconstruction method (RSIR).

Based on local conservative formulation, attempts have been done to build a HLLC-type Riemann solver. We believe that failure came from complexity of interface conditions and multiplicity of dispersed phase contact waves. Moreover, as shown hereafter, the solution is not self-similar. The same shock tube test case as before is reconsidered and the solution is shown at various times in Fig. III.3.

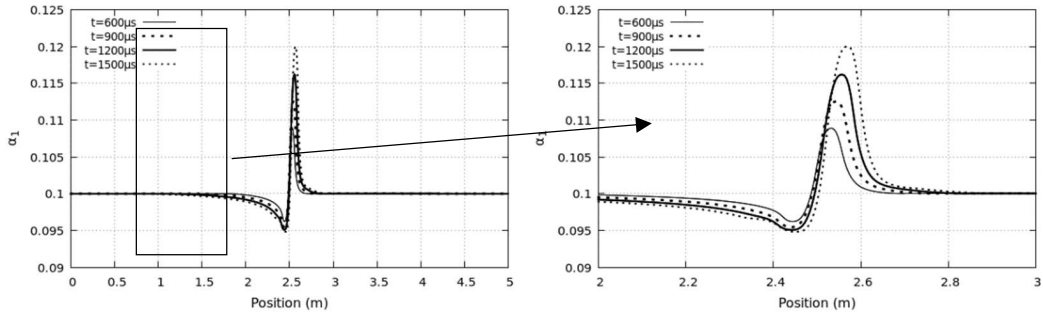

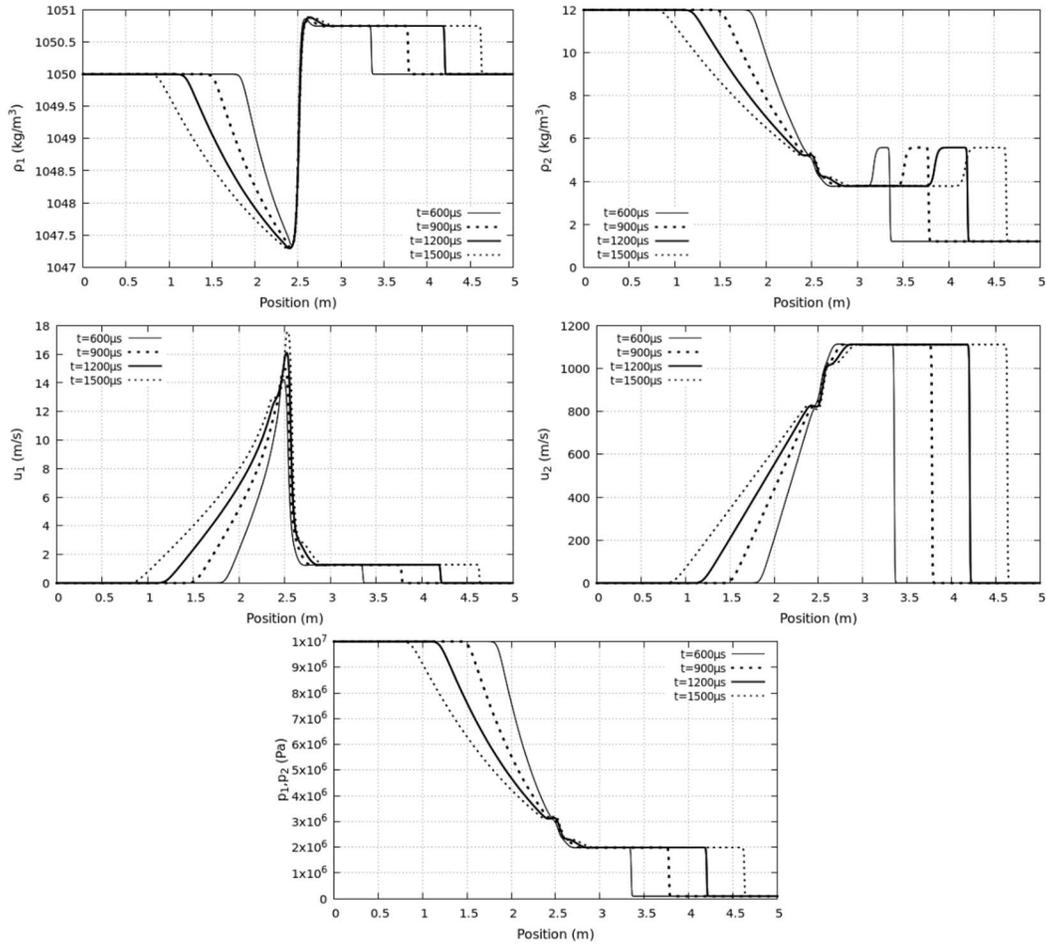

Figure III.3 – Shock tube test problem of Figure III.2 considered at various times. Length of the domain has been increased to 5 m to avoid interaction with the boundaries, and the initial discontinuity is placed at 2.5 m. Stiff pressure relaxation is used at any time. 1 000 computational cells are used with CFL=0.5. Computations are done with the MUSCL method and Minmod limiter. Solutions for the carrier phase (2) appear self-similar but solutions of the dilute phase 1 are not, regarding volume fraction and velocity. Volume fraction of the dilute phase keeps increasing without converging towards a constant state. The same tendency is observed in the velocity of the dilute phase. These observations are mesh and solver independent.

Origin of this interesting behavior is examined hereafter with the help of the Baer and Nunziato (1986) model. The same shock tube test problem is rerun with this model with the help of the HLLC-type solver of Furfaro and Saurel (2015). In the computations of Figure III.4 pressure relaxation is absent.

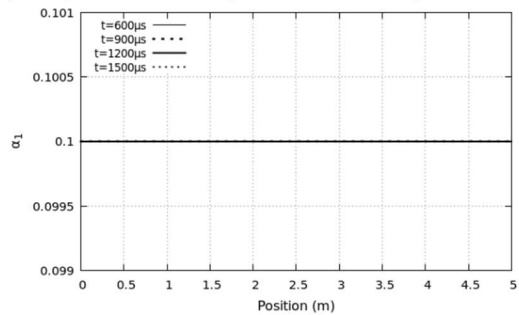

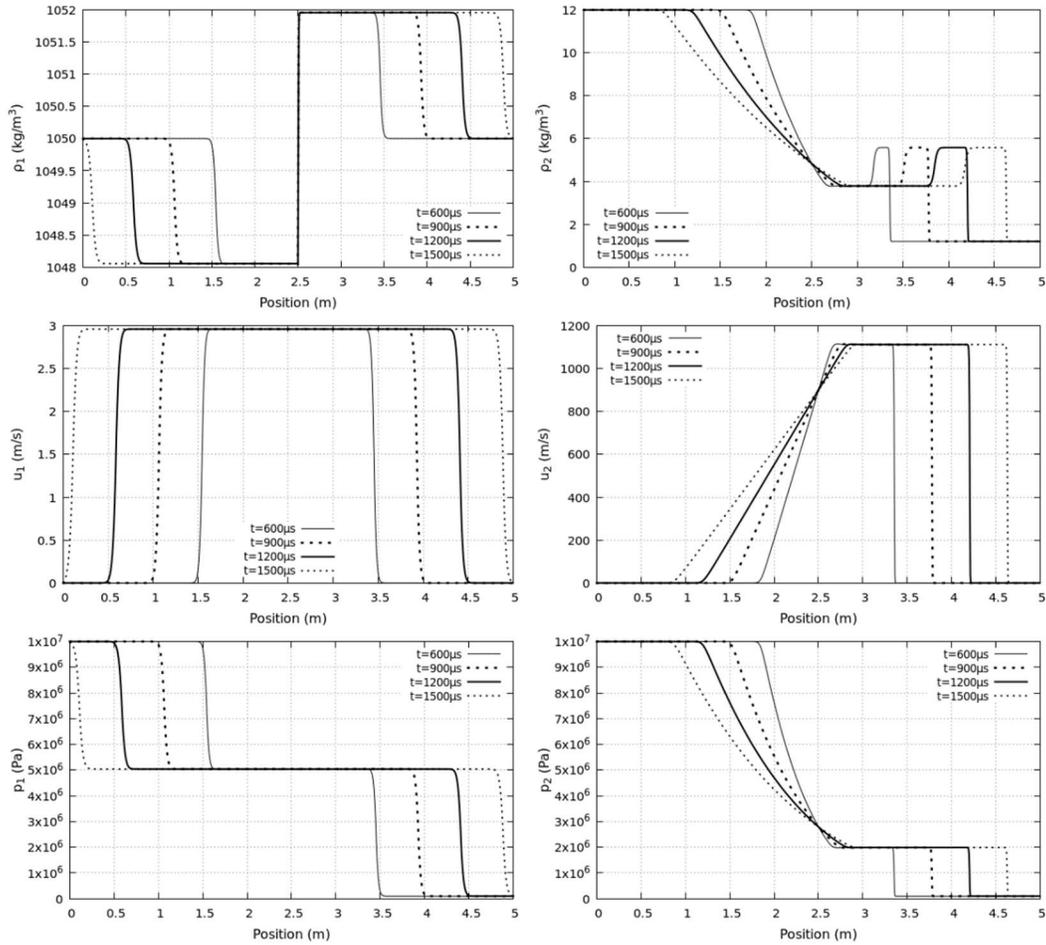

Figure III.4 – Two-phase shock tube problem of Figure III.2 computed with the Baer and Nunziato (1986) model in the absence of pressure relaxation. Length of the domain has been increased to 5m to avoid waves interaction with the boundaries, and the initial discontinuity is placed at 2.5m. 1 000 computational cells are used with CFL=0.5. Computations are done with the MUSCL method and Minmod limiter. The solution is now self-similar and consists in two decoupled shock tube solutions, as well known.

The same shock tube test problem is reconsidered once more in the presence of stiff pressure relaxation in the Baer and Nunziato (1986) model. Corresponding results are shown in Figure III.5.

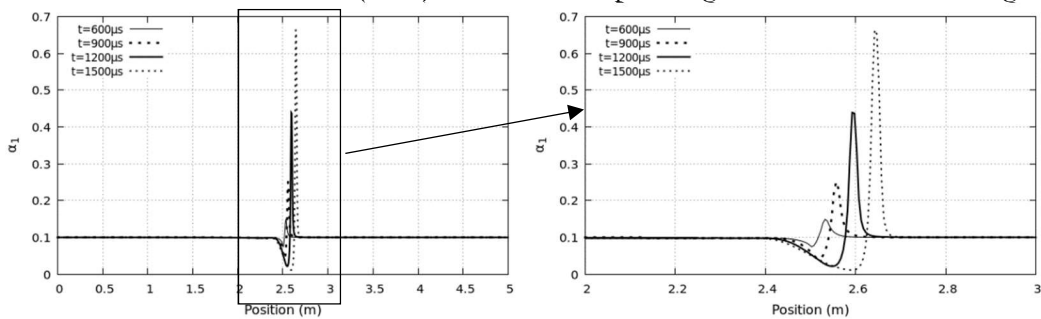

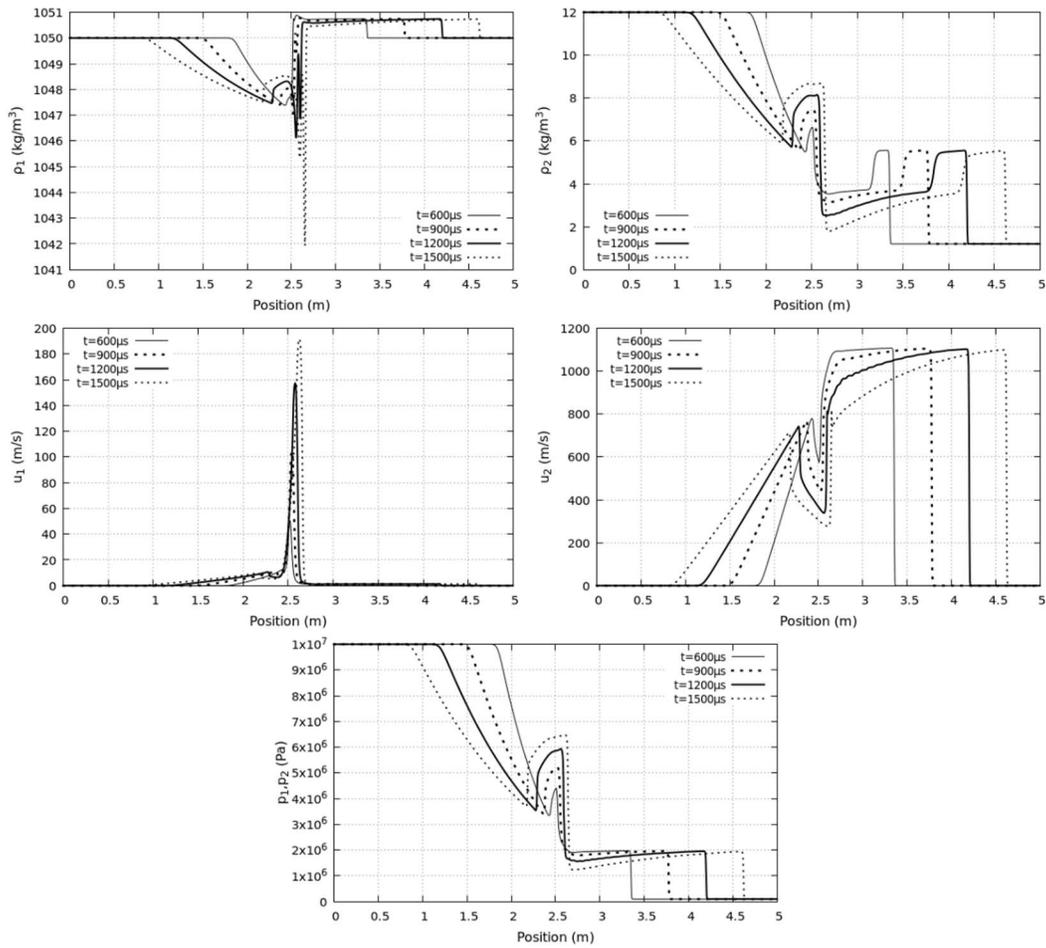

Figure III.5 – Two-phase shock tube problem of Figure III.2 computed with the Baer and Nunziato (1986) model in the presence of stiff pressure relaxation. Length of the domain has been increased to 5m to avoid waves interaction with the boundaries, and the initial discontinuity is placed at 2.5m. 1 000 computational cells are used with CFL=0.5. Computations are done with the MUSCL method and Minmod limiter. The solution is not self-similar. This is not surprising as source terms related to pressure relaxation are present.

In the presence of pressure relaxation, solutions of the BN model are no longer self-similar, and this is not surprising. It explains why solutions of flow model (III.1) are not self-similar, as this flow model has sense only in the stiff pressure relaxation limit.

Another interesting feature appears with solutions of Figures III.3 and III.5. Solution for the velocity of the dispersed phase appears multivalued, as schematized in Figure III.6. The left and right star velocities are different, implying creation of a particle's cluster, with a Dirac function type volume fraction profile. Such behavior also happens with the BN model (Figure III.5) as a combination of pressure relaxation and non-conservative terms, acting as a drag force between phases.

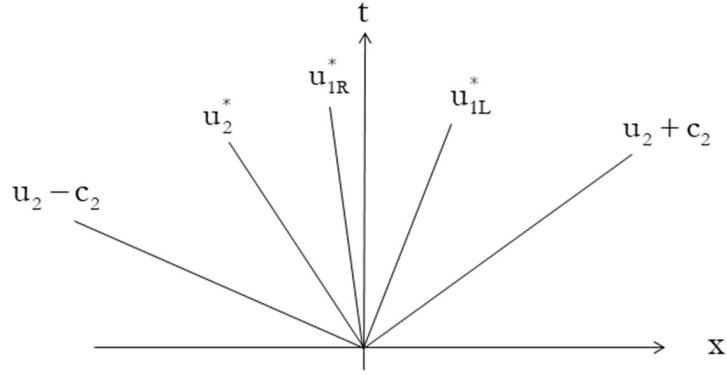

Figure III.6 – Schematic representation of the multivalued phase 1 velocity. Differences in the right and left star velocities of the dispersed phase may result in particle agglomeration as time evolve. This effect is present in the computations of Fig. III.3 and III.5.

These difficulties are omitted when dealing with single state Riemann solver, such as HLL or Rusanov. From this robust solution basis, we now address internal reconstruction to reduce numerical smearing.

### III.2 Riemann solver with internal reconstruction (RSIR) for the two-phase model

The present flow solver is not based on variations across the various waves but only on rebuilding two intermediate states to preserve isolated volume fraction discontinuities and reduce artificial smearing during transport. In this direction, the intermediate wave speed is based on the phase 1 contact wave:

$$S_{M1} = \frac{U^*_{HLL}\left((\alpha\rho u)_1\right)}{U^*_{HLL}\left((\alpha\rho)_1\right)}, \qquad (III.10)$$

with $U^*_{HLL}$ given by (III.11) in the context of System (III.6) and same wave speed estimates as before. As detailed with the Euler equations the method to solve System (III.1) proceeds in two steps:
- Determine average state with HLL based on (III.6);
- Rebuild the solution.

Thanks to the local conservative formulation the first step is immediate. System (III.6) is expressed as,

$$\frac{\partial U}{\partial t} + \frac{\partial \Phi(U)}{\partial x} = 0,$$

with $p_I$ given by (III.5).

The average state is obtained as,

$$U^*_{HLL} = \frac{\Phi_R - \Phi_L + S_L U_L - S_R U_R}{S_L - S_R}. \qquad (III.11)$$

The next step consists in the approximation of states L* and R*. Knowledge of the intermediate wave speed $S_{M1}$ from (III.10) enables following decomposition of the average state,

$$U^*_{HLL} = \omega_R U^*_R + \omega_L U^*_L \qquad (III.12)$$

with,

$$\omega_R = \frac{S_R - S_{M1}}{S_R - S_L} \quad \text{and} \quad \omega_L = \frac{S_{M1} - S_L}{S_R - S_L}.$$

As System (III.12) involves two unknown variable vectors, another set of relations is needed.

**Reconstruction**

For the dispersed phase, the Rankine-Hugoniot relations (III.3) imply, $\alpha_{1R}^* = \alpha_{1R}$ and $\alpha_{1L}^* = \alpha_{1L}$.
The difference of these two relations reads,
$$\alpha_{1R}^* - \alpha_{1L}^* = \alpha_{1R} - \alpha_{1L}.$$
Parameter $\beta$ is introduced to control numerical diffusion,
$$\alpha_{1R}^* - \alpha_{1L}^* = \beta(\alpha_{1R} - \alpha_{1L}).$$
As phase 1 density has no jump across left- and right-facing waves, similar relation is obtained:
$$(\alpha\rho)_{1R}^* - (\alpha\rho)_{1L}^* = \beta\left((\alpha\rho)_{1R} - (\alpha\rho)_{1L}\right).$$
Regarding momentum jump across the intermediate wave, the same relation is used as in the context of the Euler equations,
$$(\alpha\rho u)_{1R}^* - (\alpha\rho u)_{1L}^* = \beta\left((\alpha\rho)_{1R} - (\alpha\rho)_{1L}\right)S_{M1}.$$
The energy jump relation is based on the interface condition related to the momentum equation of phase 1:
$$(\alpha\rho)_{1R}^* u_{1R}^* \left(u_{1R}^* - S_{M1}\right) + \alpha_{1R}^* \left(p_{1R}^* - p_I\right) = (\alpha\rho)_{1L}^* u_{1L}^* \left(u_{1L}^* - S_{M1}\right) + \alpha_{1L}^* \left(p_{1L}^* - p_I\right).$$
Velocity jump conditions across left- and right-facing waves of the dispersed phase are introduced as,
$$u_{1R}^* = u_{1R} \text{ and } u_{1L}^* = u_{1L}.$$
Momentum jump across the contact wave consequently reads,
$$(\alpha\rho)_{1R}^* u_{1R} \left(u_{1R} - S_{M1}\right) + \alpha_{1R}^* \left(p_{1R}^* - p_I\right) = (\alpha\rho)_{1L}^* u_{1L} \left(u_{1L} - S_{M1}\right) + \alpha_{1L}^* \left(p_{1L}^* - p_I\right).$$
Assuming fluid 1 governed by the stiffened gas EOS, it becomes,
$$\alpha_{1R}^* \left[(\gamma_1 - 1)\rho e - \gamma_1 p_{1,\infty}\right]_{1R}^* - \alpha_{1L}^* \left[(\gamma_1 - 1)\rho e - \gamma_1 p_{1,\infty}\right]_{1R}^* = \left((\alpha\rho)_{1L}^* u_{1L}\left(u_{1L} - S_{M1}\right) - (\alpha\rho)_{1R}^* u_{1R}\left(u_{1R} - S_{M1}\right)\right)$$
$$+ \left(\alpha_{1R}^* - \alpha_{1L}^*\right)p_I$$

The internal energy jump thus reads,
$$(\alpha\rho e)_{1R}^* - (\alpha\rho e)_{1L}^* = (\alpha\rho)_{1L}^* u_{1L} \frac{(u_{1L} - S_{M1})}{(\gamma_1 - 1)} - (\alpha\rho)_{1R}^* u_{1R} \frac{(u_{1R} - S_{M1})}{(\gamma_1 - 1)} + \left(\alpha_{1R}^* - \alpha_{1L}^*\right)\frac{p_I + \gamma_1 p_{1,\infty}}{(\gamma_1 - 1)}.$$
The total energy jump for phase 1 follows,
$$(\alpha\rho E)_{1R}^* - (\alpha\rho E)_{1L}^* = \left(\alpha_{1R}^* - \alpha_{1L}^*\right)\frac{p_I + \gamma_1 p_{\infty,1}}{(\gamma_1 - 1)} + \left[(\alpha\rho)_{1R}^* - (\alpha\rho)_{1L}^*\right]\frac{S_{M1}^2}{2}$$
$$+ (\alpha\rho)_{1L}^* u_{1L} \frac{(u_{1L} - S_{M1})}{(\gamma_1 - 1)} - (\alpha\rho)_{1R}^* u_{1R} \frac{(u_{1R} - S_{M1})}{(\gamma_1 - 1)}.$$
For the first phase, these relations summarize as,
$$U_{R,k}^* = U_{L,k}^* + \psi_k, \quad k = 1,..,4, \tag{III.13}$$
with,

$$\psi = \begin{bmatrix} \beta(\alpha_{1R} - \alpha_{1L}) \\ \beta[(\alpha\rho)_{1R} - (\alpha\rho)_{1L}] \\ \beta[(\alpha\rho)_{1R} - (\alpha\rho)_{1L}]S_{M1} \\ \beta(\alpha_{1R} - \alpha_{1L})\dfrac{p_I + \gamma_1 p_{\infty,1}}{(\gamma_1 - 1)} + \beta[(\alpha\rho)_{1R} - (\alpha\rho)_{1L}]\dfrac{S_{M1}^2}{2} + (\alpha\rho)_{1L}^* u_{1L}\dfrac{(u_{1L} - S_{M1})}{(\gamma_1 - 1)} - (\alpha\rho)_{1R}^* u_{1R}\dfrac{(u_{1R} - S_{M1})}{(\gamma_1 - 1)} \end{bmatrix}$$

For the dispersed phase, combining (III.12) and (III.13) the intermediate states are computed as,

$$\begin{cases} U_{R,k}^* = U_{HLL,k}^* + \omega_L \psi_k \\ U_{L,k}^* = U_{HLL,k}^* - \omega_R \psi_k \end{cases}, \quad k = 1,..,4. \tag{III.14}$$

Reconstruction is now addressed for the second phase. For the sake of simplicity, let us consider that the volume fraction of the second phase is considered additionally. To preserve the saturation constraint ($\alpha_1 + \alpha_2 = 1$) the jump relation across the intermediate wave reads,

$$\alpha_{2R}^* - \alpha_{2L}^* = \beta(\alpha_{2R} - \alpha_{2L}).$$

Another assumption is now introduced. The density of the carrier phase is assumed uniform in the two intermediate states. In other words,

$$\overline{\rho}_2 = \dfrac{(\alpha\rho)_{2,HLL}^*}{\alpha_{2,HLL}^*}.$$

This assumption is needed as the flow solver behaves incorrectly when the product of two discontinuous functions are present (volume fraction and density). The assumption made at this level is similar to prolongated formulations used in Ghost Fluid Methods (Fedkiw et al., 1999), immersed boundary methods and diffuse interface methods (Kapila et al., 2001, Allaire et al., 2002, Massoni et al., 2002, Saurel and Pantano, 2018). In diffuse interface methods, numerical diffusion of volume fractions $\alpha_k$ and apparent densities $\alpha_k \rho_k$ results in automatic prolongment of phase's density $\rho_k$. In this frame, apparent densities are rebuilt as,

$$(\alpha\rho)_{2R}^* = \alpha_{2R}^* \overline{\rho}_2,$$
$$(\alpha\rho)_{2L}^* = \alpha_{2L}^* \overline{\rho}_2.$$

Consequently,

$$(\alpha\rho)_{2R}^* - (\alpha\rho)_{2L}^* = \beta(\alpha_{2R} - \alpha_{2L})\overline{\rho}_2.$$

Momentum jump across the intermediate wave is rebuilt as,

$$(\alpha\rho u)_{2R}^* - (\alpha\rho u)_{2L}^* = \left[(\alpha\rho)_{2R}^* - (\alpha\rho)_{2L}^*\right]S_M = \beta(\alpha_{2R} - \alpha_{2L})\overline{\rho}_2 S_{M2}$$

where the contact wave speed of the second phase is computed as,

$$S_{M2} = \dfrac{U_{HLL}^*((\alpha\rho u)_2)}{U_{HLL}^*((\alpha\rho)_2)}.$$

It remains to determine total energy jump of the second phase. At this level no distinction is made between the left and right velocities of that phase: $u_{2L}^* = u_{2R}^* = S_{M2}$. As before, the energy jump relation is based on the interface condition related to the momentum equation of phase 2:

$$(\alpha\rho)_{2R}^* S_{M2}(S_{M2} - S_{M1}) + \alpha_{2R}^*(p_{2R}^* - p_I) = (\alpha\rho)_{2L}^* S_{M2}(S_{M2} - S_{M1}) + \alpha_{2L}^*(p_{2L}^* - p_I).$$

Alternatively, it reads,

$$\alpha_{2R}^* p_{2R}^* - \alpha_{2L}^* p_{2L}^* = \left((\alpha\rho)_{2L}^* - (\alpha\rho)_{2R}^*\right) S_{M2}(S_{M2} - S_{M1}) + \left(\alpha_{2R}^* - \alpha_{2L}^*\right) p_I.$$

Assuming fluid 2 governed by the stiffened gas EOS, it becomes,

$$\alpha_{2R}^* \left[(\gamma_2 - 1)\rho e - \gamma_2 p_\infty\right]_{2R}^* - \alpha_{2L}^* \left[(\gamma_2 - 1)\rho e - \gamma_2 p_\infty\right]_{2R}^* = \left((\alpha\rho)_{2L}^* - (\alpha\rho)_{2R}^*\right) S_{M2}\left(S_{M2} - S_{M1}\right) + \left(\alpha_{2R}^* - \alpha_{2L}^*\right) p_I$$

The internal energy jump thus reads,

$$(\alpha\rho e)_{2R}^* - (\alpha\rho e)_{2L}^* = \frac{\left((\alpha\rho)_{2L}^* - (\alpha\rho)_{2R}^*\right) S_{M2}\left(S_{M2} - S_{M1}\right) + \left(\alpha_{2R}^* - \alpha_{2L}^*\right)\left(p_I + \gamma_2 p_{\infty,2}\right)}{(\gamma_2 - 1)}.$$

The total energy jump for phase 2 follows,

$$(\alpha\rho E)_{2R}^* - (\alpha\rho E)_{2L}^* = \left[(\alpha\rho)_{2R}^* - (\alpha\rho)_{2L}^*\right]\left(\frac{S_{M2}^2}{2} - S_{M2}\frac{(S_{M2} - S_{M1})}{(\gamma_2 - 1)}\right) + \left(\alpha_{2R}^* - \alpha_{2L}^*\right)\frac{(p_I + \gamma_2 p_{\infty,2})}{(\gamma_2 - 1)}.$$

The carrier phase is thus rebuilt as,

$$U_{R,k}^* = U_{L,k}^* + \psi_k, \quad k = 5,..,8, \tag{III.15}$$

with,

$$\psi = \beta(\alpha_{2R} - \alpha_{2L})\begin{bmatrix} 1 \\ \dfrac{1}{\rho_2} \\ \overline{\rho_2} S_{M2} \\ \overline{\rho_2}\left(\dfrac{S_{M2}^2}{2} - S_{M2}\dfrac{(S_{M2} - S_{M1})}{(\gamma_2 - 1)}\right) + \dfrac{(p_I + \gamma_2 p_{\infty,2})}{(\gamma_2 - 1)} \end{bmatrix}$$

Last, relations (III.14) and (III.15) are used to compute $U_{R,k}^*$ and $U_{L,k}^*$ for $k = 5,..,8$.

At this level, states vectors $U_R^*$ and $U_L^*$ are determined.

Local conservative fluxes are computed as,

$$\begin{cases} \Phi_R^* = \Phi_R + S_R\left(U_R^* - U_R\right) \\ \Phi_L^* = \Phi_L + S_L\left(U_L^* - U_L\right) \end{cases}. \tag{III.16}$$

Fluxes of System (III.1) are computed as,

$$F_k^* = \Phi_k^* + \Lambda_k \quad k = L, R, \tag{III.17}$$

with $\Lambda_k = p_I\left(0,\ 0,\ \alpha_{1k}^*,\ (\alpha u)_{1k}^*,\ 0,\ (1 - \alpha_{1k}^*),\ -(\alpha u)_{1k}^*\right)^T$.

$\alpha_k^*$ are determined form HLL state as,

$$\alpha_k^* = \frac{\Phi(\alpha_{k,R}) - \Phi(\alpha_{k,L}) + S_L \alpha_{k,L} - S_R \alpha_{k,R}}{S_L - S_R}. \tag{III.18}$$

The associated Godunov-type scheme including non-conservative terms reads,

$$U_i^{n+1} = U_i^n - \frac{\Delta t}{\Delta x}\left(F_{i+\frac{1}{2}}^* - F_{i-\frac{1}{2}}^*\right) + \Delta t H_i,$$

with $H_i$ given by,

$$H_{i,u} = p_i^n \frac{\alpha^*_{1,i+\frac{1}{2}} - \alpha^*_{1,i-\frac{1}{2}}}{\Delta x}$$

$$H_{i,E} = p_i^n \frac{\Phi^*_{i+\frac{1}{2}}(\alpha_1) - \Phi^*_{i-\frac{1}{2}}(\alpha_1)}{\Delta x}.$$

In these expressions $\alpha^*_{k,i\pm\frac{1}{2}}$ and $\Phi^*_{i+\frac{1}{2}}(\alpha_1)$ are given by (III.18) and (III.16) respectively.

### III.3 Examples and validations

Validations of the flow solver and comparisons with the former Rusanov method are addressed first. Second computational examples are shown showing method's capabilities.

a) Validation of the RSIR solver

A volume fraction discontinuity at rest is considered to check method capability to maintain such stationary wave. Corresponding results are shown in Figure III.7.

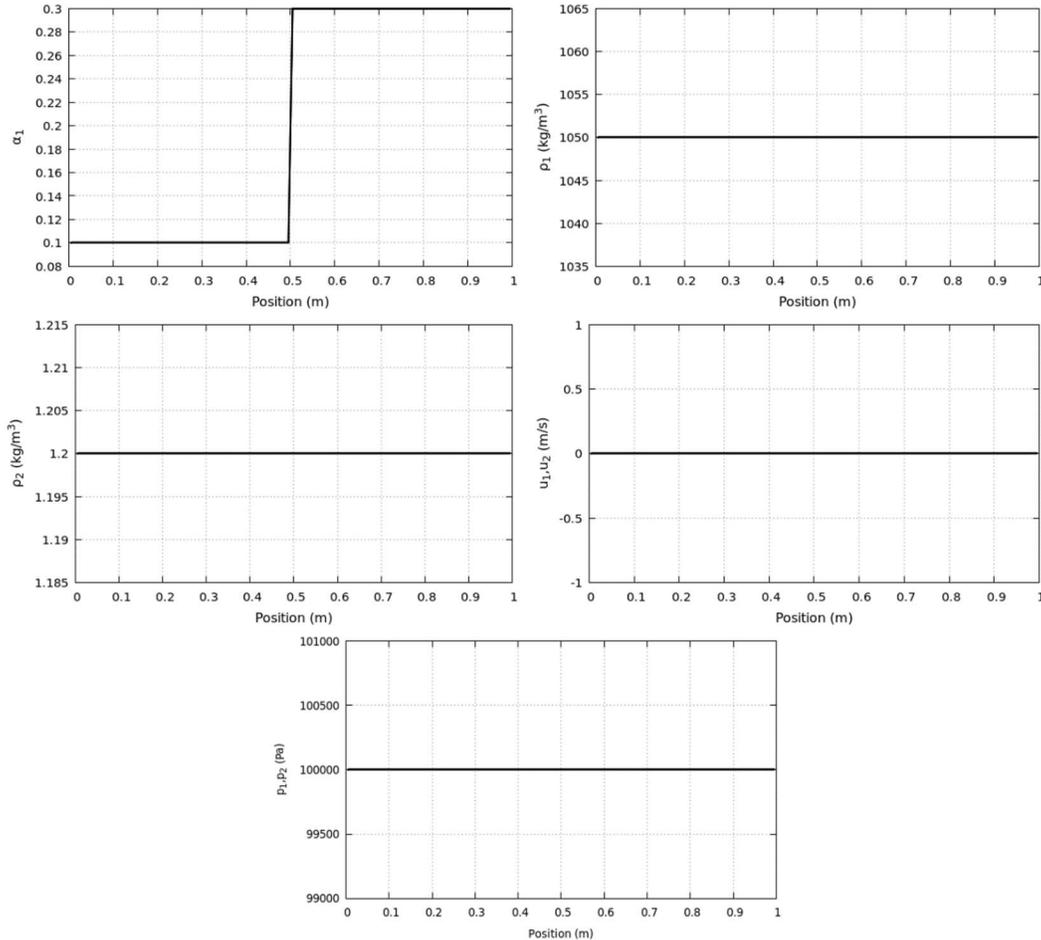

Figure III.7 – Results obtained by the new solver (β=1) embedded in the Godunov-MUSCL method with Minmod limiter for the computation of a contact discontinuity at rest. 100 computational cells are used with CFL=0.5. Results are shown at time t=6ms. The discontinuity is well preserved and spurious pressure and velocity oscillations are absent.

Transport of the same volume fraction discontinuity is now considered in a flow in uniform pressure and velocity conditions. Corresponding results are shown in Figure III.8.

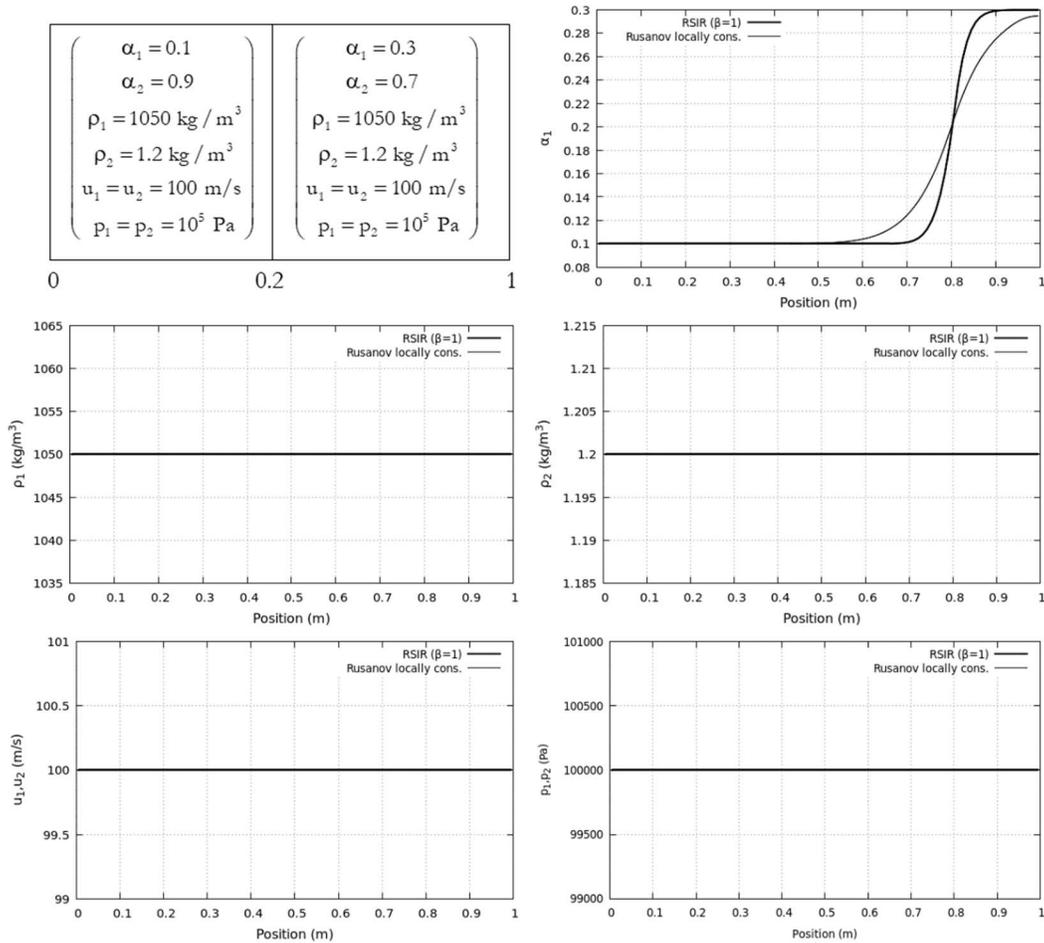

Figure III.8 – Results obtained by the new solver ($\beta=1$) embedded in the Godunov-MUSCL method with Minmod limiter for the computation of volume fraction discontinuity transport in a uniform pressure and velocity flow. 100 computational cells are used with CFL=0.5. Results are shown at time t=6ms.

The same two-phase shock tube test problem as in Figures III.2 and III.3, in the absence of drag force is now considered. The new method and the Rusanov solver are compared in Figure III.9.

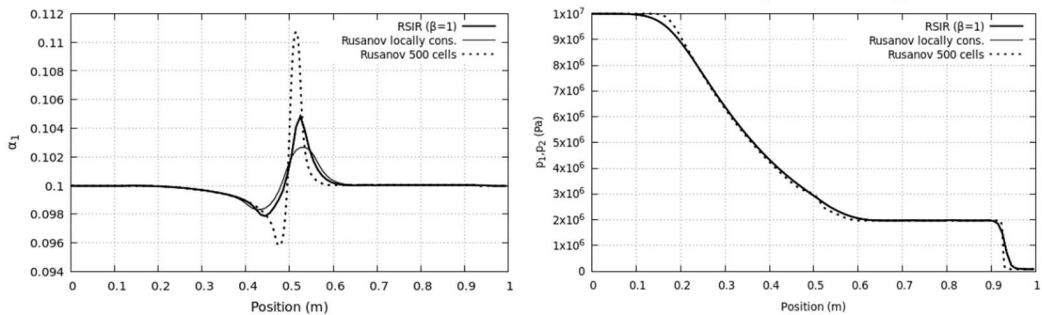

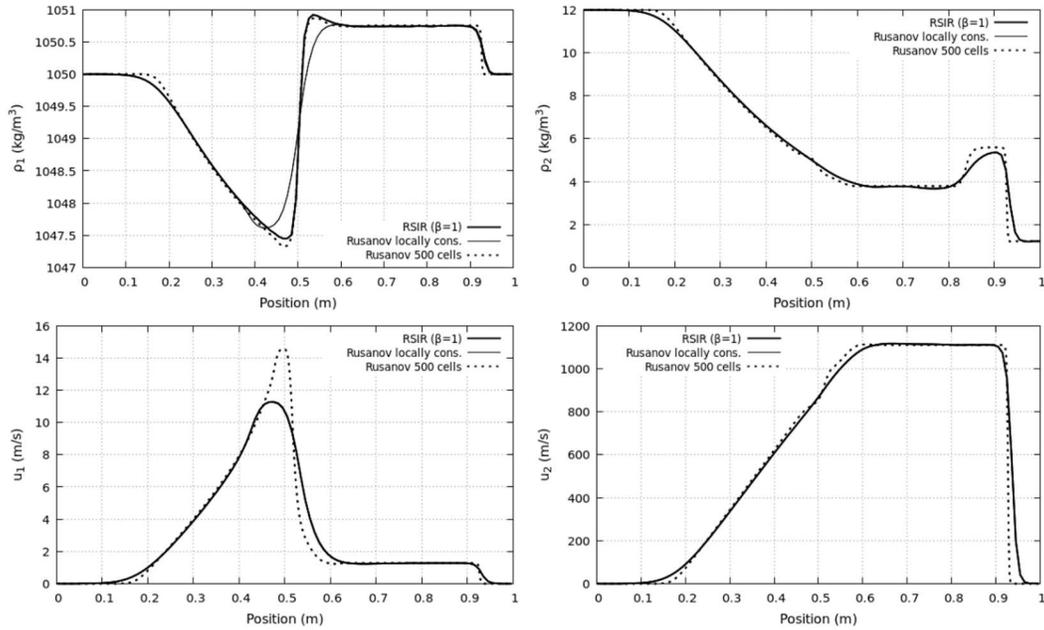

Figure III.9 – Comparison of the results obtained by the new solver (β=1) and the local conservative solver of Rusanov, both embedded in the Godunov-MUSCL method with Minmod limiter for the computation of the shock tube presented in Figure III.2. 100 computational cells are used with CFL=0.5. Results are shown at time t=300μs. Stiff pressure relaxation is used. With the RSIR method, main improvements appear in the volume fraction and phase 1 density computations. Indeed, the RSIR solution lies between computed results with the Rusanov method with 100 and 500 cells.

The RSIR solver is consequently validated and improves accuracy of the Rusanov and HLL solvers. Also, it preserves volume fraction discontinuities at rest. Its capabilities are now illustrated on the computation of a challenging two-phase flow instability.

## IV. Multi-D example: Particle jetting during radial explosion

When a spherical or cylindrical explosive charge is surrounded by a liquid layer or a granular particle bed the material dispersal occurs through particle jets having well defined size. On the example shown in the Figure IV.1 a cylindrical explosive charge is initially surrounded by a liquid layer.

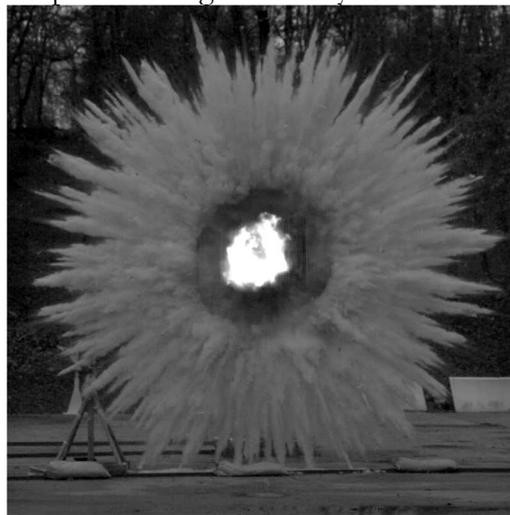

Figure IV.1 – A cylindrical explosive charge is initially surrounded by a liquid layer. When the charge explodes the liquid layer transforms to a cloud of droplets forming highly dynamical particle jets.

Gas expansion during explosion fragments the liquid layer to a cloud of droplets that form highly dynamical particle jets. The same observation is reported when a granular bed is used instead of a liquid layer. Dispersion is consequently not spherical or cylindrical in the sense that one-dimensional computations result in significant errors in the predictions of materials presence. Experimental and numerical studies of this phenomenon have been achieved by Zhang et al. (2001), Milne et al. (2010), Frost (2010), Parrish and Worland (2010) to cite a few. Simplified situations have been considered in Rodriguez et al. (2013) and Xue et al. (2018). The explosive is replaced by a shock tube and the matter to disperse is placed between two plates, in a Hele-Shaw cell. Other simplified situations have been considered for example in McGrath et al. (2018), Osnes et al. (2018), Carmouze et al. (2018) to study possible clustering effects due to aerodynamic forces. It seems that jets formation and size selection mechanism is still unidentified.

In the following a configuration like the one considered in Rodriguez et al. (2013) with a Hele-Shaw cell and a particle ring is studied. Typical results reported in Rodriguez thesis are shown in Figure IV.2 at times 5 ms, 8 ms and 57 ms after rupture of the shock tube diaphragm, inducing shock wave and gas flow through a ring of flour particles.

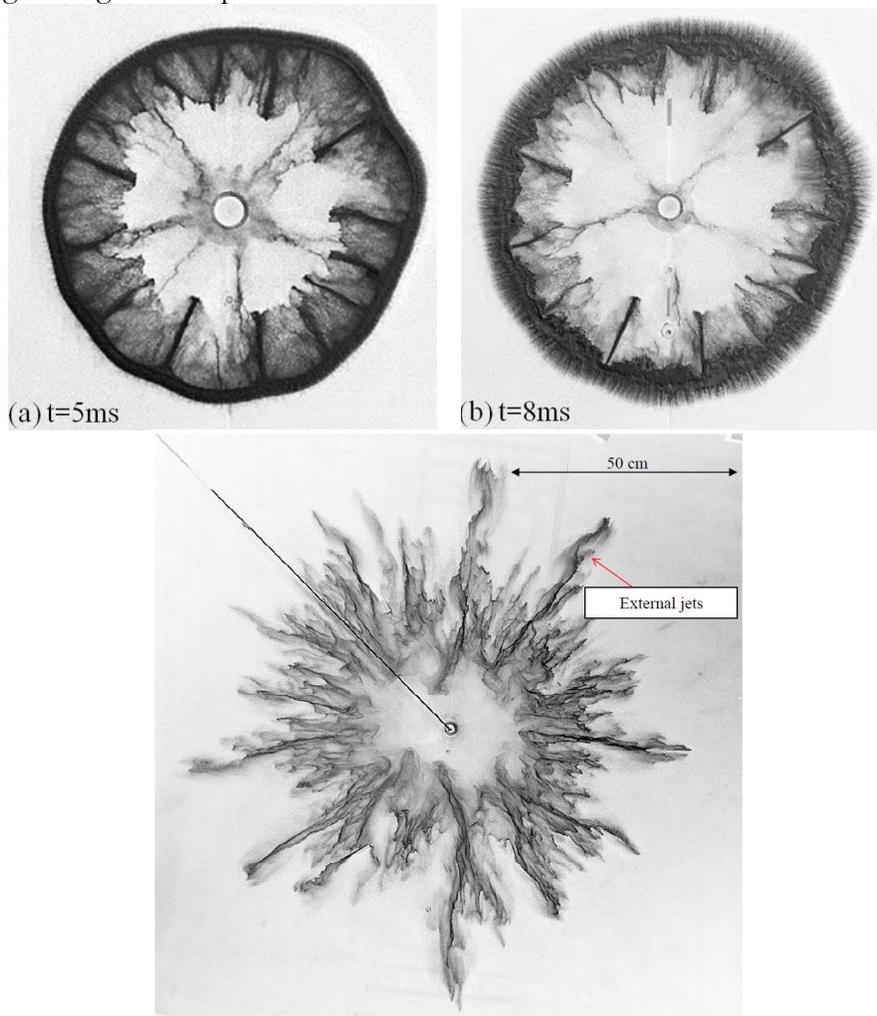

Figure IV.2 – Typical interfacial instabilities reported in Rodriguez thesis and papers, as well as in Xue et al. (2018). Impulsive motion of a particle ring by a gas flow induces well defined particles fingers flowing to the center direction, oppositely to the gas flow. At later times, here at 8 ms, short wavelength instabilities also appear at the external surface. As time evolves, external surface instabilities grow and become dominant, as shown in the third picture at time 57 ms. Internal jets are thus observed at early times, followed by external ones at late times.

As reported by Rodriguez et al. (2013) and Xue et al. (2018), instabilities appear first at the inner interface and second at the outer one. Shape of these fingers is singular, in the sense that they do not qualitatively compare to the Richtmyer-Meshov instabilities or Rayleigh-Taylor ones, nor any other known instability. Indeed, mushroom type shape is observed with these two instabilities, while fingers are observed in the present context

With the help of the new model and present RIRS solver attempt is done to reproduce at least qualitatively these instabilities. Former attempts by the authors were based on Marble (1963) model and fail to reproduce any, inner or outer, instability. In the present attempt computations are based on flow model III.1 extended to 2D and resolved numerically in the DALPHADT code on triangular cells. Compared to the Marble model, the present one has a fundamental difference. Non-conservative terms are present in the momentum and energy equations. These non-conservative terms are often called 'nozzeling terms' in reference to the Euler equations with variable cross section. We prefer howether to interpret them as 'differential drag force'. As shown for example in Figure III.9, phase 1 is set to motion as a result of these terms, in the absence of conventional drag force. Indeed, drag parameter $\lambda$ is set to zero in these computations. As shown later, it seems that non-conservative terms are responsible for appearance of these instabilities and size selection.

Considered computational domain and initial data are reported in Figure IV.3.

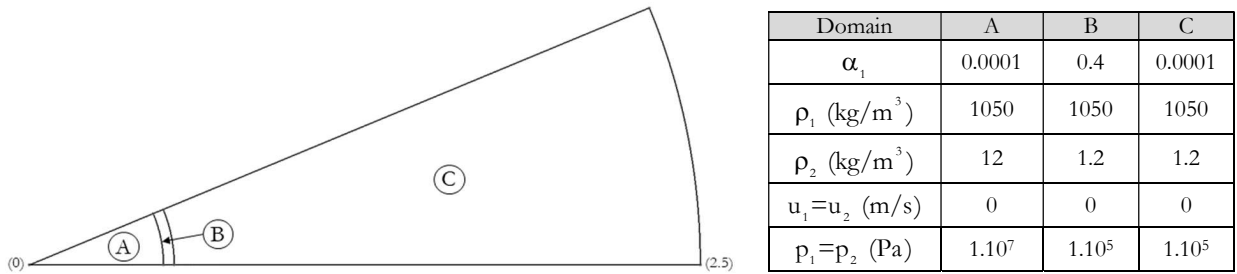

| Domain | A | B | C |
|---|---|---|---|
| $\alpha_1$ | 0.0001 | 0.4 | 0.0001 |
| $\rho_1$ (kg/m$^3$) | 1050 | 1050 | 1050 |
| $\rho_2$ (kg/m$^3$) | 12 | 1.2 | 1.2 |
| $u_1 = u_2$ (m/s) | 0 | 0 | 0 |
| $p_1 = p_2$ (Pa) | $1.10^7$ | $1.10^5$ | $1.10^5$ |

Figure IV.3 – Computational domain and initial data. The portion represents 1/16 of the complete disc with angle $\theta = \pi/8$. The domain denoted A corresponds to the high-pressure chamber, filled with gas only. The domain denoted B represents the initial particles ring of 4cm width. The domain denoted C corresponds to the low-pressure chamber at atmospheric conditions.

The mesh is made of triangles and contains about 15 cells in the radial direction of the particles ring.

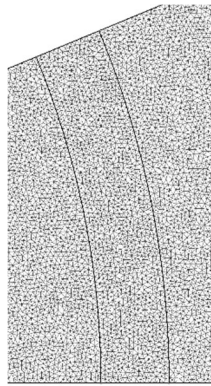

Figure IV.4 – Representation of the mesh size used in the particle jetting computations. In the entire domain 361 222 computational cells are used.

Typical computed results are shown in Figures IV.5 and IV.6 and IV.7 at different times. Volume fraction contours are shown. These computations show that the flow model reproduces, at least qualitatively both fingering instabilities issued of the inner and outer surfaces. In these computations drag force is considered with constant particles diameter (1 mm) and air kinematic viscosity (18.10$^{-6}$ Pa.s). Drag force is modelled through Clift and Gauvin (1971) correlation,

$$F_D = \frac{3}{8R_1}\alpha_1 C_d \rho_2 \|\vec{u}_2 - \vec{u}_1\|(u_2 - u_1),$$

with,

$$C_d = \begin{cases} \dfrac{24}{Re_1}\left(1 + 0.15 Re_1^{0.687}\right) & \text{if } Re_1 < 800 \\ 0.438 & \text{otherwise} \end{cases}$$

and $Re_1 = \dfrac{2R_1 \rho_2 \|\vec{u}_2 - \vec{u}_1\|}{\mu_2}$.

Drag force $F_D$ is inserted in the right-hand side of the momentum equation of phase 1 and its opposite in phase 2. The power of this force $F_D.u_1$ is inserted similarly in associated energy equations.

In the absence of drag effects internal particles jets selection is observed as well but quantitative differences appear at later times. It is worth to mention that particle volume fraction is initially uniform in the cloud. No arbitrary wavelength is introduced in the initial data. Computed results at early times are shown in Figure IV.5.

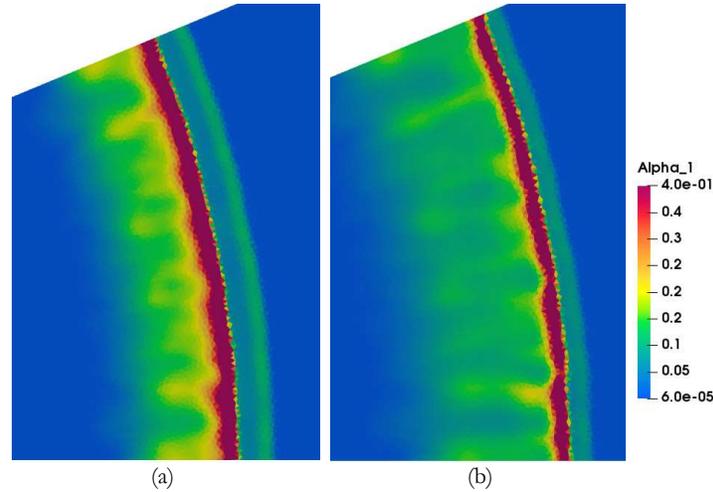

(a)        (b)

Figure IV.5 – Volume fraction contours of the dispersed phase for the particle jetting simulation, focused on the particles cloud at early times: (a) t=0.6 ms and (b) t=0.75 ms. 361 222 cells are used with CFL=0.9. Results are obtained with the RSIR (β=1) solver embedded in the MUSCL scheme with Superbee limiter. A compaction zone appears first in the cloud in red color. Particles jets develop at the inner interface and direct to the center domain. Their growth is visible by comparing their length in graphs (a) and (b). They qualitatively look like the instabilities observed in Figure IV.2 (a) and (b). Another front appears at the outer surface but appears more like a diffusion zone rather than the short wavelength instabilities visible in Figure IV.2 (b). This is possibly due to insufficient special resolution of the present computations.

In Figure IV.5 internal jets appear and develop to the center domain direction. A compaction zone appears in the cloud and a detached front with low particles concentration also appears ahead the compaction front. Noticeably the compaction front in red color catches up the detached front. Evolution at intermediate times is reported in Figure IV.6.

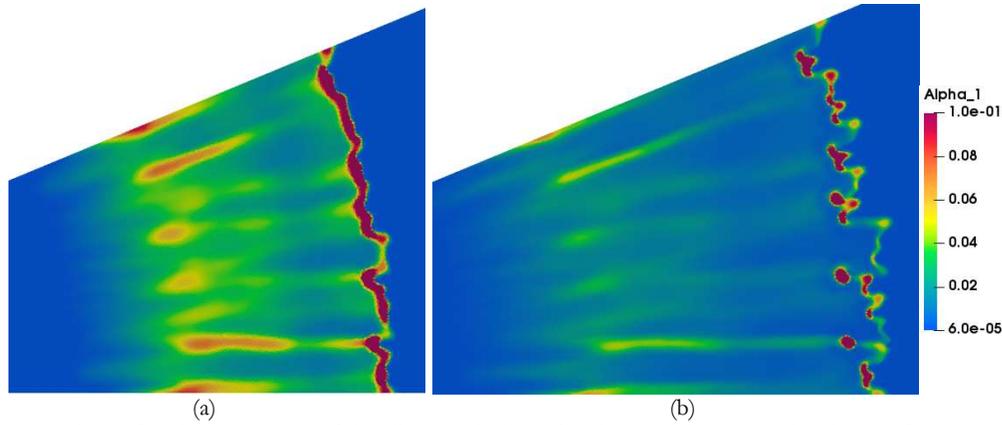

Figure IV.6 – Volume fraction contours of the dispersed phase for the particle jetting simulation, focused on the cloud at intermediate times: (a) t=1.5 ms and (b) t=2.25 ms. Same computational parameters as those of Figure IV.5 are used as it corresponds to the same computation. The compaction front and the detached one are now merged and destabilize. Particles concentration zones having cluster type shapes appear in graph (b). Inner jets are still present and continue development.

In Figure IV.6 the external front destabilizes, and relatively dense particle clusters appear. The inner front jets flowing to the domain center continue their development. Evolution at later times is reported in Figure IV.7.

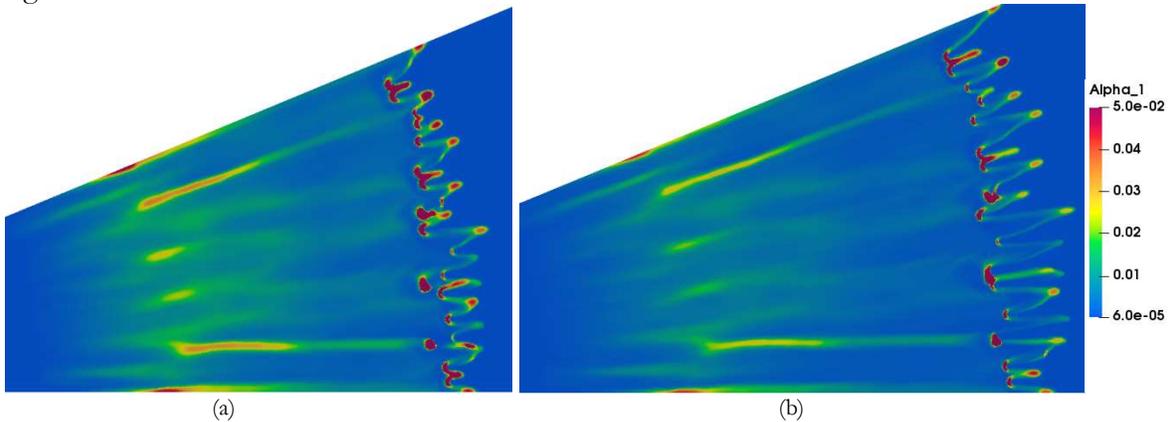

Figure IV.7 – Volume fraction contours of the dispersed phase for the particle jetting simulation, focused on the particles cloud at later times: (a) t=3 ms and (b) t=3.75 ms. Same computational parameters as those of Figure IV.5 and IV.6 are used. External front instabilities are now created and develop. Dilution of the internal jets happens while external jets develop as a consequence of particles 'dense' zones created at intermediate times. External jets amplitude grows as visible by comparing results of graphs (a) and (b).

In Figure IV.7 external front instabilities are created and develop while internal ones tend to vanish. It is also interesting to note that the number of external fingers varies between graphs (a) and (b) as reported by Rodriguez et al. (2013) regarding experimental observations. In graph (b) about 10 leading fingers emerge while in (a) they are about 14.

Although not precisely identified from the present numerical experiments, the formation mechanism of this fingering instability appears closely related to non-conservative terms. They play the role of a differential drag force, acting intensively at cloud boundaries and vanishing in the wake, when volume fraction gradients disappear.

At the modeling level, non-conservative terms present similar form as capillary ones (Brackbill et al., 1992, Perigaud and Saurel, 2005) except that curvature effects are absent in the present two-phase formulation. Another major difference is that cloud boundaries are obviously highly permeable in the present context, while interfaces are not permeable in conventional hydrodynamic instabilities, except those considering flames and phase transition, where low permeability is present.

## V. Conclusion

A Riemann solver with internal reconstruction (RSIR) has been built as an extension of Linde (2002) solver. It has been first developed for the Euler equations and shown to provide similar accuracy and robustness as the HLLC solver, while being not as systematic as this one. The method has been secondly extended to a two-phase non-equilibrium model developed by the authors. This model presented serious difficulties as it is weakly hyperbolic and valid only in the presence of stiff pressure relaxation, rendering solutions not self-similar. Thanks to the RSIR approach, a low dissipation solver has been developed. It has been validated against solution obtained with more conventional, but dissipative solvers. The new method has been applied in the last section to a difficult problem of fingering instability in granular media and has shown possible explanation of the formation mechanism. Extra work is still needed in this special two-phase flow topic to achieve understanding of this instability.

Regarding the Riemann solver, it seems flexible for many applications where most of the physics is governed by the two extreme waves and an intermediate one. Moreover, a parameter is present to control dissipation when flow conditions are particularly severe. In all present computations the solver appeared robust with parameter ($\beta=1$), corresponding to the minimum dissipation.

**Acknowledgements.** Part of this work has been carried out in the framework of the Labex MEC (ANR-10-LABX-0092) and of the A*MIDEX project (ANR-11-IDEX-0001-02), funded by the Investissements d'Avenir French Government program managed by the French National Research Agency (ANR). We also acknowledge funding from ANR through project SUBSUPERJET ANR-14-CE22-0014.